\documentclass{article}
\usepackage[utf8]{inputenc}
\usepackage{mathrsfs}
\usepackage{authblk}
\usepackage{mathtools}
\usepackage[colorinlistoftodos]{todonotes}  

\usepackage{comment}
\usepackage{amsmath}
\usepackage{amssymb}
\usepackage{amsthm}
\numberwithin{equation}{section}

\usepackage[left=3cm, right=3cm]{geometry}

\usepackage[T1]{fontenc}
\usepackage[utf8]{inputenc}
\usepackage[english]{babel}
  \usepackage{lmodern}  
\usepackage{graphicx}
\usepackage{subfigure}
\usepackage{tikz}
\usepackage{makeidx}
\makeindex

\usepackage{multicol}
\usepackage{float}
\usepackage{tikz,pgf}
\usepackage{tikz-cd}
\usepackage{float}
\usepackage{xcolor}
\usepackage[all,knot,arc,color,web]{xy}
\usepackage{graphicx}
\xyoption{curve}
\usepackage{bbold}
\usepackage{fancyhdr}
\usepackage{stmaryrd}
\usepackage{hyperref}
\theoremstyle{plain}
\newtheorem{theorem}{Theorem}[section]
\newtheorem*{theorem*}{Theorem}
\newtheorem{proposition}[theorem]{Proposition}
\newtheorem*{proposition*}{Proposition}
\newtheorem{corollary}[theorem]{Corollary}
\newtheorem{lemma}[theorem]{Lemma}
\newtheorem{definition}[theorem]{Definition}

\theoremstyle{definition}
\newtheorem{remark}[theorem]{Remark}

\newtheorem{example}[theorem]{Example}

\title{Birkhoff attractors for dissipative symplectic billiards}
\author[1]{Luca Baracco}
\author[1]{Olga Bernardi}
\author[2]{Anna Florio}
\author[1]{Alessandra Nardi}
\affil[1]{Dipartimento di Matematica Tullio Levi--Civita, Università di Padova,
via Trieste 63, 35121 Padova, Italy.}
\affil[2]{\textsc{Ceremade}-Universit\'e Paris Dauphine-PSL\\
75775 Paris, France.}
\date{\today}

\begin{document}

\maketitle

\abstract{\noindent The aim of the present paper is to propose and study a dissipative variant of symplectic billiards within planar strictly convex domains. The associated billiard map is dissipative, thus it admits a compact invariant set, the so-called Birkhoff attractor. Its complexity depends on the rate of the dissipation as well as on the geometry of the billiard table. We prove that $(a)$ for strong dissipation, the Birkhoff attractor is a normally contracted graph over the zero section; $(b)$ for mild dissipation, the Birkhoff attractor within a centrally symmetric domain is an indecomposable continuum whose restricted dynamics has positive topological entropy. We compare these results with the case of dissipative Birkhoff billiards, studied in \cite{BFL}.}
\section{Introduction}
\noindent Symplectic billiards were introduced by P. Albers and S. Tabachnikov in 2018 \cite{AT} and were subsequently studied by various authors who investigated integrability \cite{BaBe}, \cite{BBN2}, its Mather's $\beta$-function \cite{BBN1} and area spectral rigidity \cite{BBN3}, \cite{FSV}, \cite{Tsod}.

Let $\Omega \subset \mathbb{R}^2$ be a strictly convex domain with $C^k$ boundary $\partial \Omega$, $k \ge 2$. Fix an origin $O \in \mathrm{int}(\Omega)$ and an orientation of $\partial\Omega$. Let $\mathbb{S}=\mathbb{R}/2\pi\mathbb{Z} \ni t \mapsto \gamma(t) \in \mathbb{R}^2$ be a $C^k$ parametrization of $\partial \Omega$ such that $\gamma(\mathbb{S}) = \partial \Omega$. Given $t_1,t_2\in\mathbb{S}$, if $\gamma(t_1)$ and $\gamma(t_2)$ are two successive bounces, then the next bounce under the symplectic dynamics occurs at $\gamma(\tilde{t}_3)$ if and only if the vector $\gamma(\tilde{t}_3)-\gamma(t_1)$ is parallel to $\gamma'(t_2)$.

The associated billiard map $T$ is a twist map, which preserves an area form and whose generating function is
$$L(t_1,t_2) := \det(\gamma(t_1) - O, \gamma(t_2) - O)\, .$$
We refer to Section \ref{INTRO DEF} for all details. In the present paper, we introduce a definition of \textit{dissipative} symplectic billiard map, as explained right below.

\begin{figure}[H]
 \centering
\includegraphics[scale=0.35]{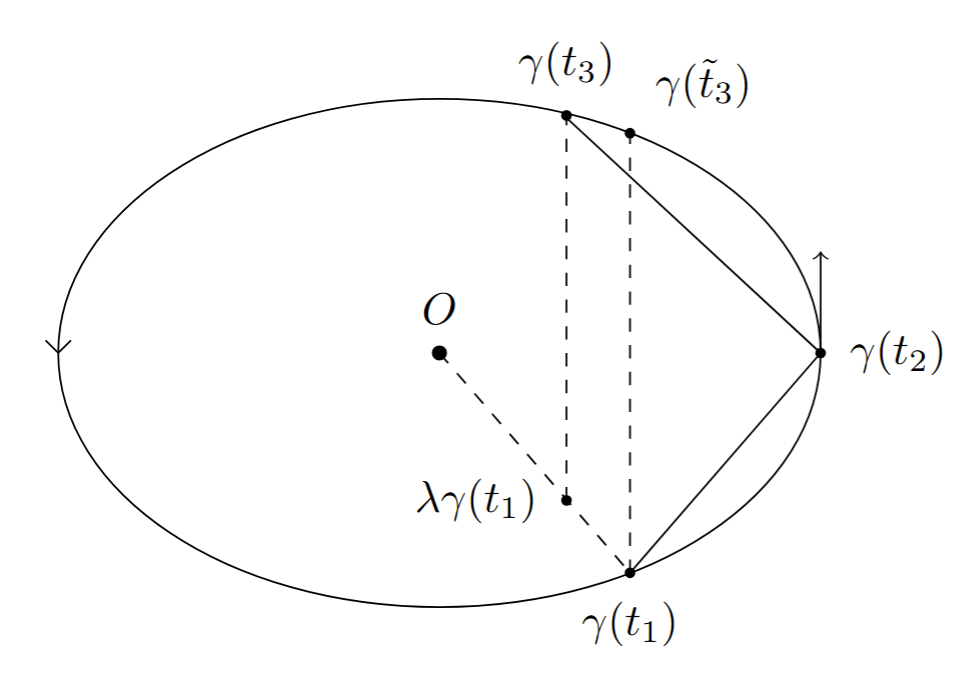}
    \caption{The dissipative symplectic billiard map compared with the conservative one.}
    \label{DSB}
\end{figure}

Fix $\lambda \in (0,1]$. If $\gamma(t_1)$ and $\gamma(t_2)$ are two successive bounces, then the dissipative symplectic dynamics gives $\gamma(t_3)$ as next bounce if and only if the vector $\gamma(t_3) - \lambda \gamma(t_1)$ is parallel to $\gamma'(t_2)$. 
In other words, the segment passing through $\gamma(t_3)$ and parallel to $\gamma'(t_2)$ is approaching the origin, when $\lambda$ approaches $O$. See Figure \ref{DSB}. We will see that the dissipative billiard map $T_\lambda$ is still a twist map, but it no longer preserves an area form. Clearly $T_1 = T$. 

\begin{figure}[H]
 \centering
\includegraphics[scale=0.4]{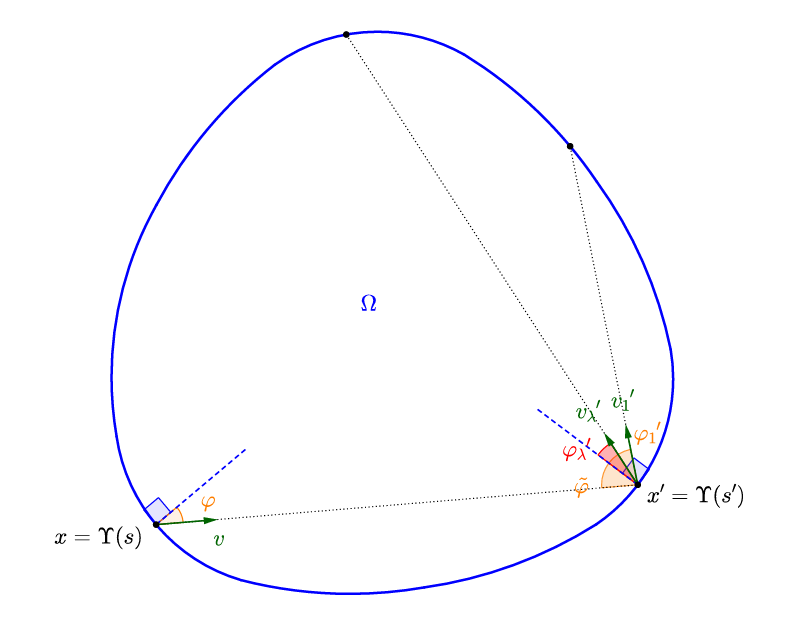}
    \caption{The dissipative Birkhoff billiard map compared with the conservative one.}
    \label{DBB}
\end{figure}

\noindent A dissipative version of standard 
(Birkhoff) billiards has been recently introduced by Bernardi, Florio and Leguil in \cite{BFL}. In such a case, the usual reflection law, which requires that the angle of incidence equals the angle of reflection, is changed in such a way that the reflected angle bends toward the inner normal at the incidence point. See Figure \ref{DBB}. Moreover, different notions of billiard with some form of dissipation have previously been considered by other authors, see for example \cite{Day}, \cite{Tab}, \cite{MPS}, $\dots$ 

Let us now move into some properties of dissipative symplectic billiards, see Section \ref{INIZIO} for all details. We first remark that --up to an appropriate choice of coordinates-- the phase space for the symplectic billiard map $T$ is a bounded cylinder, denoted by $\mathcal{P}$. In the coordinates $(t,s)\in\mathcal{P}$, it results that 
$$T_\lambda(t,s) = \mathcal{H}_\lambda \circ T(t,s)\, ,$$ 
where $\mathcal{H}_\lambda$ is the $\lambda$-contraction map along the fiber, i.e., $\mathcal{H}_\lambda(t,s) = (t,\lambda s)$. Clearly, the map $T_\lambda$ is no longer conservative; in fact, it is dissipative in the sense of Le Calvez \cite{LeCal}, see Definition \ref{DEF DISS}. In particular, $T_\lambda$ contracts the standard area form, that is $T^*_\lambda (dt \wedge ds) = \lambda dt \wedge ds$. According to the vocabulary of symplectic dynamics, the map $T_\lambda$ is conformally symplectic. Due to its dissipative character, the map $T_\lambda$ admits a \textit{global attractor} ``\`a la Conley''
$$\Lambda_0 := \bigcap_{n\in\mathbb{N}}T_\lambda^n(\mathcal{P})\, ,$$
which is not empty, compact, connected and $T_\lambda$-invariant. Moreover, it separates the phase space $\mathcal{P}$, that is its complement $\mathcal{P}\setminus \Lambda_0$ is a disjoint union of two connected open sets $U_-$ and $U_+$. Following Birkhoff \cite{Birkhoff} and Le Calvez \cite[Page 91]{LeCalEns}-- it is possible to detect the smallest compact, connected, $T_\lambda$-invariant set that separates $\mathcal{P}$: the so-called \textit{Birkhoff attractor}. Obtained from $\Lambda_0$ ``by removing the whiskers'', it is then defined as
$$\Lambda := \mathrm{cl}(U_-) \cap \mathrm{cl}(U_+)\, .$$
We need to underline that, differently from $\Lambda_0$ and despite its name, in general $\Lambda$ is not an attractor in the usual sense, that is it is not the omega-limit set of one of its neighborhoods. We refer to Section \ref{SUB 2.2} for all details on the Birkhoff attractor $\Lambda$. 

The study of the properties and the structure of the Birkhoff attractor for dissipative maps of the annulus, begun by Birkhoff \cite{Birkhoff}, pursued with Charpentier \cite{Char34}, Le Calvez \cite{LeCal,LeCalEns}, Crovisier \cite{Cro}, \dots\ Different criteria to obtain chaotic invariant continua for dissipative maps are proved in \cite{BG,Koropecki,Alejandro,PassTal,PassTalPlus}, \dots\ Recently, the notion of Birkhoff attractor has been generalized in higher dimensions by Arnaud, Humili\`ere and Viterbo, see \cite{AHV,Vit22,VitPlusPlus}, for conformally symplectic maps on cotangent bundles. This generalization relies on symplectic invariants of exact Lagrangian submanifolds.

More generally, in the last years, conformally symplectic dynamics has been an active research area, see \cite{AF,AFR,AA}, \dots\ The dissipative systems associated to the discounted Hamilton-Jacobi equation have been studied in \cite{Zavi,Zavi2,Zavi3,AHV}, showing interesting properties of solutions of such equations through Weak KAM methods. In the dissipative Hamiltonian framework, invariant KAM-like tori for conformally symplectic flows are studied by Calleja, Celletti and de la Llave \cite{CCDLL}; very recently, Gidea, de la Llave and Seara \cite{GDLLS} have studied the geometry and the structure of normally hyperbolic manifolds and their scattering maps for conformally symplectic Hamiltonian flows.\\

The aim of the present paper is to study the dynamical and topological complexity of the Birkhoff attractor in terms of the rate of the dissipation $\lambda \in (0,1)$ as well as in terms of the geometry of the billiard table.

In Section \ref{INIZIO} we start by discussing some domains whose corresponding dissipative symplectic billiard dynamics $T_\lambda$ has --independently on $\lambda \in (0,1)$-- a particularly simple Birkhoff attractor, that is $\Lambda = \mathbb{S} \times \{0\}$. This is the case of
any (centrally symmetric) Radon domain $\Omega$, once fixed the origin $O$ to be the center of symmetry, see Proposition \ref{Radon doamins} and Remark \ref{RADON and ELLIPSE}. It is worth noting that this is a first substantial difference compared to dissipative Birkhoff billiards, whose Birkhoff attractor is assured to be the zero section only in the case of the circle. \\
\indent Section \ref{NORMALLY} is devoted to study $\Lambda$ for dissipative symplectic billiards when the dissipation is strong, i.e., $0<\lambda\ll 1$. In such a case, independently on the choice of the origin, $T_{\lambda}$ exhibits a topologically simple Birkhoff attractor, that is $\Lambda$ coincides with $\Lambda_0$ and it is a normally contracted graph over $\mathbb{S}$, as precised in the following statement. 

\begin{theorem*}Let $\Omega\subset \mathbb{R}^2$ be a strongly convex (i.e., with never vanishing curvature) domain with $C^k$ boundary, $k \ge 2$. 
 \begin{itemize}
\item[$(a)$] There exists $\lambda(\Omega)\in(0,1)$ such that, for $\lambda\in(0,\lambda(\Omega))$, the Birkhoff attractor $\Lambda$ coincides with the global attractor $\Lambda_0$ and it is a normally contracted $C^1$ graph over $\mathbb{S}$. 
\item[ $(b)$] There exists $\lambda'(\Omega)<\lambda(\Omega)$ such that, for $\lambda\in (0, \lambda'(\Omega))$, the Birkhoff attractor $\Lambda$ is a $C^{k-1}$ graph over $\mathbb{S}$ and it converges to $\mathbb{S}\times\lbrace0\rbrace$, as $\lambda\to0$, in the $C^1$ topology.
\end{itemize}
\end{theorem*}

We refer to Theorem \ref{nc} for all details. The proof is based on a cone-field criterion (see Proposition \ref{cones}) and standard results in normally hyperbolic dynamics, and it follows the same lines of the proof of the corresponding result for dissipative Birkhoff billiards (see \cite[Theorem 5.7]{BFL}). However, in the Birkhoff dissipative dynamics, the analogous theorem does not hold in general, but only for a class of strictly convex domains satisfying a geometric \textit{pinching} condition (see \cite[Definition D]{BFL}), while in the symplectic case, no hypothesis on the geometry of the table is needed. \\
\indent In Section \ref{CSSB}, we focus on dissipative symplectic billiard maps for a centrally symmetric domain $\Omega$. In particular, we prove that, for a rate of dissipation small enough (i.e., when the dissipation is strong), not only the Birkhoff attractor is a graph, but we can describe quite explicitely the dynamics restricted to it. We refer to Propositions \ref{punti in sezione} and \ref{unione varieta instabili} for the detailed statements and we resume here the main result.

\begin{proposition*}Let $\Omega$ be a strongly convex centrally symmetric domain with $C^k$ boundary, $k \ge 2$. 
\begin{itemize}
\item[$(a)$] There exists $\lambda(\Omega)\in (0,1)$ such that, for $\lambda \in (0,\lambda(\Omega))$, the Birkhoff attractor $\Lambda$ intersects the zero section exactly in the 4-periodic points. 
\item[$(b)$] There exists $\lambda'(\Omega)\in (0,1)$ such that, for $\lambda \in (0,\lambda'(\Omega))$, there exists an open and dense set of centrally symmetric domains, whose associated Birkhoff attractor has rotation number $1/4$ and decomposes as 
\[\Lambda =\bigcup_{i=1}^l \bigcup_{j=0}^3\overline{\mathcal{W}^u(T_\lambda^j(H_i);T_\lambda^4)}\, ,\]
 where $\lbrace H_i\rbrace_{i=1}^l$ is a finite family of $4$-periodic points of saddle type, and $\mathcal{W}^u(H;T_\lambda^4)$ is the unstable manifold of a saddle point $H$, with respect to the dynamics of $T_\lambda^4$.
\end{itemize}
\end{proposition*}

For dissipative Birkhoff billiards, the corresponding results are Corollary 3.4 and Theorem 5.14 in \cite{BFL}. Even if the main ideas of the proofs in the two settings have common points, for dissipative symplectic billiards there is no need --as noticed above-- of any geometric pinching condition of the table. Moreover, in the case of dissipative Birkhoff billiards, 
axial symmetry (instead of central symmetry) and $2$-periodic orbits (instead of $4$-periodic orbits) play a fundamental role. 

\indent Section \ref{birkhoff-section-complicated} focuses on the Birkhoff attractor under conditions of weak dissipation, i.e., when the rate of dissipation is close to $1$. In contrast with the previous cases, we obtain examples of complex Birkhoff attractors. 
A sufficient condition in order to observe topologically and dynamically intriguing phenomena is that $T_1 = T$ admits an instability region containing the zero section $\mathbb{S} \times \{0\}$. This follows by an adaptation of a result of Le Calvez (here Proposition \ref{lecalvez}): in such a case, the Birkhoff attractor $\Lambda$ for the corresponding dissipative dynamics, if $\lambda$ is close enough to 1, admits different upper and lower rotation numbers $\rho^- < \rho^+$ , dynamical quantities defined in Section \ref{SUB 2.2}. This property has various consequences for $\Lambda$, all already observed in \cite{BFL} for dissipative Birkhoff billiards. In particular: 
\begin{itemize}
\item[$(i)$] $\Lambda$ is an \textit{indecomposable continuum}, that is it cannot be written as the union of two compact, connected, non trivial sets (from a result by Charpentier, see \cite{Char34}). 
\item[$(ii)$] For every rational $\frac p q \in (\rho^- , \rho^+)$, there exists a periodic point in $\Lambda$ of rotation number $\frac p q$ (which can be deduced from \cite{BG}).
\item[$(iii)$] If $x$ is a saddle periodic point of rotation number $\frac p q$, for $\frac p q\in(\rho^-,\rho^+)$, then its unstable manifold is contained in $\Lambda$ (see \cite[Proposition 14.3]{LeCal}).
\item[$(iv)$] The map $T_{\lambda}$ restricted to $\Lambda$ has positive topological entropy, as a consequence of the existence of a rotational horseshoe (see \cite[Theorem A]{Alejandro}).
\end{itemize}
~\newline
\noindent The original part of the section is therefore devoted to prove when hypotheses of Proposition \ref{lecalvez} --hence guaranteeing ``complicated'' Birkhoff attractors-- are satisfied, as resumed here below. We refer to Theorem \ref{N PER} and Proposition \ref{mather docet} respectively. 

\begin{theorem*}Let $\Omega\subset \mathbb{R}^2$ be a strictly convex domain with $C^k$ boundary, $k\geq 2$. The conservative symplectic billiard map $T$ associated to $\Omega$ admits an instability region containing the zero section $\mathbb{S} \times \{0\}$ in the following two cases:
\begin{itemize} 
\item[$(a)$] If $\Omega$ belongs to an open and dense set of strongly convex, centrally symmetric billiard tables.
\item[$(b)$] If $\Omega$ has at least one point of zero curvature. 
\end{itemize}
\end{theorem*}

Case $(a)$ is based on these two facts. Among strongly convex, centrally symmetric tables, an essential invariant curve for $T$ passing through the zero section has necessarily rotation number $1/4$, see Proposition \ref{Hera}; however, such curves are very easy to destroy with an arbitrary small perturbation of the table, see Proposition \ref{open and dense no 1 4 inv curve}. We stress that, regarding the case $(a)$, the hypothesis of central symmetry plays a fundamental role. This is another difference from dissipative Birkhoff billiards where, for weak dissipation, no properties on the geometry of the billiard table are required, see \cite[Proposition 6.15]{BFL}. 
Case $(b)$, as for Birkhoff billiards, is a straightforward consequence of Mather's theorem on the non existence of caustics.

In Section \ref{NUM SIM}, we present some numerical simulations to illustrate the above results. By using \textit{Mathematica}, we compute the billiard map $T_\lambda$ for specific domains, both centrally symmetric and not, and plot some orbits in the corresponding phase space.

Finally, we would like to emphasize that, beyond the applications to Birkhoff attractors, the paper presents and proves many properties for the map $T_{\lambda}$ (with $\lambda \in (0,1]$, so also for the conservative case), including: the formula for the differential of $T_\lambda$ (Lemma \ref{formula DT lambda}),
the discussion of the quality and quantity of $4$-periodic orbits for $T$ in a centrally
symmetric domain (Lemma \ref{bifurcation} and Lemma \ref{lemma1}), the result on the fragility of invariant curves of rotation number $1/4$ for $T$ (see Proposition \ref{open and dense no 1 4 inv curve}), \dots\ We consider that these results may be useful for further studies on symplectic billiard dynamics.
\vspace{10pt}

\noindent \textit{Acknowledgements.} A.F. is partially supported by the ANR project CoSyDy (ANR-CE40-0014), the ANR project GALS (ANR-23-CE40-0001) and PEPS project ``Jeuns chercheurs et jeunes chercheuses'' 2025. L.B. and O.B. have been partially supported by the PRIN project 2022FPZEES ``Stability in Hamiltonian Dynamics and Beyond''.

\section{Preliminaries}

This section is devoting to introduce two types of dynamics, namely, the dynamics of symplectic billiards, and the dynamics of dissipative maps on the $2$-dimensional annulus. The aim of the present work is defining a combination of these two dynamics, i.e., \emph{dissipative symplectic billiards}, and studying the main properties of the associated (Birkhoff) attractor. The model of dissipative symplectic billiard map is proposed in Subsection \ref{INIZIO}.

\subsection{(Birkhoff) attractors for dissipative maps} \label{SUB 2.2}

In this subsection, we briefly recall the definition of attractor and Birkhoff attractor for a dissipative map of the annulus. We refer to \cite{LeCal} for further details. \\
\noindent Denote $\mathbb{S}:=\mathbb{R}/2\pi \mathbb{Z}$ and $\pi \colon \mathbb{R} \to \mathbb{S}$ a universal covering. Given $\varphi_{\pm}: \mathbb{S}\to \mathbb{R}$ two continuous functions such that $\varphi_-<\varphi_+$, let $\mathscr{C}\subset \mathbb{S}\times\mathbb{R}$ be the open, relatively compact subset 
\[
\mathscr{C}:=\{ (t,s)\in \mathbb{S}\times\mathbb{R} :\ \varphi_-(s)<t<\varphi_+(s) \}\, .
\]
Consider the area form $dt\wedge ds$ on the annulus $\mathbb{S}\times\mathbb{R}$ and let $m$ be the Lebesgue induced measure. Given a set $A$, we denote by $\mathrm{int}(A)$ the interior of the set, and by $\mathrm{cl}(A)$ the closure of the set. 

\begin{definition} \label{DEF DISS}
A map $f\colon \mathscr{C}\to \mathrm{int}(\mathscr{C})$ is a dissipative map if:
\begin{enumerate}
    \item $f$ is an homeomorphism of $\mathscr{C}$ into its image, homotopic to the identity;
    \item $f$ is a $C^1$ diffeomorphism of $\mathrm{int}(\mathscr{C})$ into its image;
    \item there exists $\lambda\in (0,1)$ such that, for every $(t,s)\in\mathrm{int}(\mathscr{C})$, we have $0< \mathrm{det}Df(t,s)\leq\lambda$.
\end{enumerate}
\end{definition}
\begin{example}\label{es-semplice}
Clearly, the easiest example of dissipative map one can think of is 
$$f\colon \mathbb{S}\times [-1,1] \ni (t,s) \to (t,\lambda s) \in \times [-1,1]\, ,$$
where $\lambda\in(0,1)$.
\end{example}

\noindent Given a dissipative map $f$, since $f(\mathscr{C})\subset \mathrm{int}(\mathscr{C})$, the following definition is well-posed.

\begin{definition} 
Let $f\colon \mathscr{C}\to \mathrm{int}(\mathscr{C})$ be a dissipative map. The \emph{attractor} of $f$ is the set
\[
\Lambda_0=\Lambda_0(f):= \bigcap_{n\in\mathbb{N}}f^n(\mathscr{C})\, .
\]
\end{definition}
\noindent Observe that $\Lambda_0$ is a compact, non-empty, connected, $f$-invariant set. Moreover, it \emph{separates} the annulus $\mathscr{C}$, i.e., its complementary set $\mathscr{C}\setminus\Lambda_0$ is the union of two open, disjoint set $U_+\sqcup U_-$, such that the graph of $\varphi_\pm$ is contained in $\partial U_\pm$. The definition of the attractor depends on the initial domain $\mathscr{C}$. Moreover, the dynamics restricted to the attractor could be a priori further decomposed into smaller invariant pieces, i.e. the attractor could be not ``minimal''. This intuition --which will be clarified in Proposition \ref{Birk-attr-minimal}--  justifies the following definition, introduced by Birkhoff in \cite{Birkhoff}.
\begin{definition}
 Let $f\colon \mathscr{C}\to \mathrm{int}(\mathscr{C})$ be a dissipative map and let $\Lambda_0$ be its attractor. In particular, we have $\mathscr{C}\setminus \Lambda_0=U_-\sqcup U_+$. The \emph{Birkhoff attractor} of $f$ is then
 \[
 \Lambda=\Lambda(f):= \mathrm{cl}(U_-)\cap \mathrm{cl}(U_+)\, .
 \]
\end{definition}
\begin{proposition}\label{Birk-attr-minimal}
Let $f\colon \mathscr{C}\to \mathrm{int}(\mathscr{C})$ be a dissipative map. A set $A\subset \mathscr{C}$ separates the annulus if $\mathscr{C}\setminus A$ is the disjoint union of two open sets, containing the upper and lower boundary of $\mathscr{C}$ respectively. Denote by $\mathcal{X}$ the set of compact, non-empty, $f$-invariant subsets of $\mathscr{C}$ that separate the annulus. Then, the Birkhoff attractor $\Lambda\in \mathcal{X}$ and it is the smallest one with respect to the inclusion.
\end{proposition}
\begin{figure}[ht]
    \centering
    \includegraphics[scale=0.7]{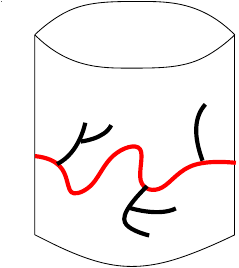}
    \caption{The red bold line corresponds to the Birkhoff attractor, while the black and red is the attractor.}
    \label{Birk-no-attr}
\end{figure}

\begin{remark}
Clearly, for the dissipative map of Example \ref{es-semplice}, the attractor coincides with the Birkhoff attractor, i.e. $\Lambda_0=\Lambda= \mathbb{S}\times\{0\}$. However, it is important to observe that the Birkhoff attractor is not a priori an attractor. There could be examples when the Birkhoff attractor is strictly contained in the attractor, as shown in Figure \ref{Birk-no-attr}. 
\end{remark}

\noindent So far, the definition of (Birkhoff) attractor does not need any twist property: the intuition is that we obtain an interesting dynamical set, which a priori should have ``lower" dimension (actually it has zero measure). Nevertheless, it could still be quite complicated, as we will see in Section \ref{birkhoff-section-complicated}. The detection of its possible topological complexity passes through the notion of upper and lower rotation number of the Birkhoff attractor. In order to have useful definitions, the map $f$ needs to be a dissipative twist map, so we continue by recalling the definition of (exact) twist map. \\
~\newline
Denote by $p_1\colon \mathbb{S}\times\mathbb{R}\to \mathbb{S}$ and $p_2\colon \mathbb{S}\times\mathbb{R}\to \mathbb{R}$ the projections over the first and second coordinates, respectively. 
Endow $\mathbb{S}\times\mathbb{R}$ with the 2-form $dt\wedge ds$: then, the 1-form $-s\, dt$ is a primitive of $dt\wedge ds$. 
\begin{definition} \label{aggiunta anna}
Let $f\colon \mathbb{S}\times\mathbb{R}\to \mathbb{S}\times\mathbb{R}$ be a $C^1$ diffeomorphism, homotopic to the identity. We say that $f$ is a \emph{positive} (resp. \emph{negative}) \emph{twist map} if there exists $\epsilon>0$ such that, for every $(t,s)\in\mathbb{S}\times\mathbb{R}$, we have
    \[
\dfrac{\partial p_1\circ f}{\partial s}(t,s) >\epsilon \qquad (\text{resp. $<-\epsilon$})\, .
    \]
Moreover, a twist map $f$ is called conservative if $f^*(dt\wedge ds )=dt\wedge ds$ and $f$ is exact if the 1-form $f^*(-s\, dt)+s\, dt$ is exact.
\end{definition}

In the sequel, we will largely use the notion of \emph{generating function} for a conservative twist map. We introduce it here through the following proposition (see \cite[Proposition 1.8]{MCA}).

\begin{proposition}\label{prop gen func}
 Let $F\colon \mathbb{R}^2\to \mathbb{R}^2$ be a $C^1$ diffeomorphism. Then, $F$ is a lift of a conservative twist map $f$ if and only if there exists a $C^2$ function $S\colon \mathbb{R}^2\to \mathbb{R}$ such that
 \begin{enumerate}
     \item for every $T,\tilde T\in\mathbb{R}$, one has $S(T+2\pi, \tilde T+2\pi) = S(T,\tilde T)$;
     \item there exists $\epsilon>0$ such that for all $T,\tilde T\in \mathbb{R}$ one has $\epsilon < -\frac{\partial^2 S}{\partial T\partial\tilde T}(T,\tilde T)$;
     \item we have $F(T,s)=(\tilde T,\tilde s)$ $\Leftrightarrow$ $\tilde s = \frac{\partial  S}{\partial \tilde T}(T,\tilde T)$ and $s=-\frac{\partial S}{\partial T}(T,\tilde T)$.
 \end{enumerate}
 We say that $S$ is a generating function for $F$ or $f$.
\end{proposition}

The notion of twist map and of generating function can be adapted also to the case of a map on a bounded annulus, as it will be the case for symplectic (conservative) billiards.

We return now to Birkhoff attractors. Since $\Lambda$ separates the annulus, the set $\mathscr{C}\setminus \Lambda$ can be written as $U^\Lambda_-\sqcup U^\Lambda_+$, where $U^\Lambda_\pm$ is an open set containing the graph of $\varphi_\pm$ in its boundary. For every $(t,s)\in\mathscr{C}$, the upper and lower vertical lines are
\[
V^+(t,s):=\{(t,y)\in\mathscr{C} :\ y\geq s\} \quad\text{and}\quad V^-(t,s):=\{(t,y)\in\mathscr{C} :\ y\leq s\} \, .
\]
Define now
\[
\Lambda^+:=\{(s,t)\in \Lambda :\ V^+(t,s)\setminus\{(t,s)\}\subset U_+^\Lambda\}
\]
and
\[
\Lambda^-:=\{(s,t)\in \Lambda :\ V^-(t,s)\setminus\{(t,s)\}\subset U_-^\Lambda\}\, .
\]
Recalling that $\pi\colon \mathbb{R}\to \mathbb{S}$ is a universal covering of $\mathbb{S}$, we use $\Pi\colon \tilde{\mathscr{C}}\subset \mathbb{R}^2\to \mathscr{C}\subset \mathbb{S}\times\mathbb{R}$ for the induced universal covering on the considered annulus. With an abuse of notation, we denote by $p_1, p_2$ the projections on the first and second coordinate, respectively, on both $\mathbb{S}\times\mathbb{R}$ and $\mathbb{R}^2$. Let $F$ be a continuous lift of $f$, where $f$ is a dissipative, positive twist map.

\begin{proposition}
    Let $f\colon \mathscr{C}\to \mathrm{int}(\mathscr{C})$ be a dissipative, positive twist map. Let $\Lambda$ be its Birkhoff attractor. The sequence
    \[
    \left(\dfrac{p_1\circ F -p_1}{n}\right)_{n\in\mathbb{N}}
    \]
    converges uniformly on $\Pi^{-1}(\Lambda^+)$ (resp. on $\Pi^{-1}(\Lambda^-)$) to a constant $\rho^+$ (resp. $\rho^-$). This constant is called upper (resp. lower) rotation number.
\end{proposition}
\noindent The next result is due to Charpentier \cite{Char34}: it provides a sufficient condition on the upper and lower rotation numbers for the existence of a “complicated” Birkhoff attractor.
\begin{theorem}
Let $f\colon \mathscr{C}\to \mathrm{int}(\mathscr{C})$ be a dissipative map. If $\rho^+-\rho^->0$, then the corresponding Birkhoff attractor is an indecomposable continuum, 
i.e., it cannot be written as the union of two compact, connected, non-trivial sets.
\end{theorem}

\subsection{Symplectic billiards} \label{INTRO DEF}
\noindent Let $\Omega\subset \mathbb{R}^2$ be a strictly convex planar domain with $C^k$ boundary $\partial \Omega$, $k\geq 2$. Assume that 
the perimeter of $\partial \Omega$ is normalized to $2 \pi$. 
Fix the origin $O \in \mathrm{int}(\Omega)$ and the positive counter-clockwise orientation on $\partial\Omega$. Let 
$$\mathbb{S} \ni t \mapsto \gamma(t) \in \mathbb{R}^2$$ 
be a $C^k$ parametrization of $\partial \Omega$, such that $\gamma(\mathbb{S})=\partial\Omega$. Given two vectors $v_1,v_2\in\mathbb{R}^2$, we denote by $\det(v_1,v_2)$ the determinant of the matrix whose columns are $v_1,v_2$: it corresponds to the signed area of the parallelogram determined by vectors $v_1$ and $v_2$. With an abuse of notation, given a point $x\in \mathbb{R}^2$, we will think of it also as a vector, considering $x-O$. Thus, given two points $x_1,x_2\in\mathbb{R}^2$, the notation $\det(x_1,x_2)$ corresponds to $\det(x_1-O,x_2-O)$.
\begin{figure}[ht]
 \centering
\includegraphics[scale=0.35]{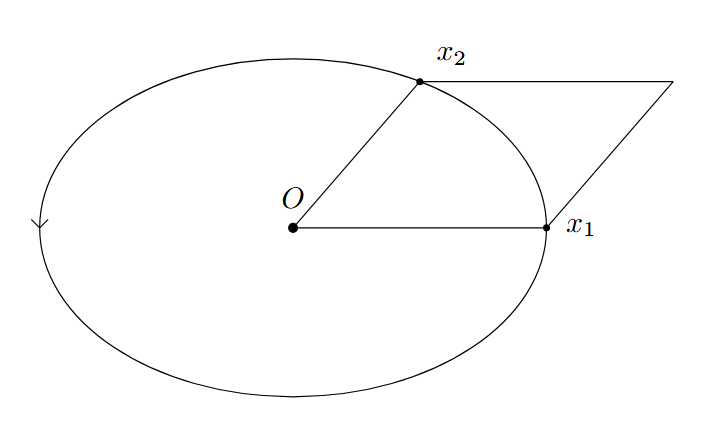}
    \caption{The determinant $\det(x_1-O,x_2-O)$ is the area of the parallelogram in figure.}
    \label{L}
\end{figure}

Let us indicate then
\[
L\colon (v,u)\in\mathbb{R}^2\to L(v,u):=\det(v,u)\in\mathbb{R}\, .
\]
In particular, for all $t_1,t_2\in\mathbb{S}$, the notation \[
L(\gamma(t_1),\gamma(t_2))
\]
denotes the signed area of the parallelogram of sides $\gamma(t_1)-O$ and $\gamma(t_2)-O$.

\noindent Given a point $\gamma(t)\in\partial\Omega$, we denote by $\gamma'(t)$ the tangent vector to $\partial\Omega$ at the point $\gamma(t)$, with respect to the fixed parametrization.
Since $\Omega$ is  strictly convex, for every point $\gamma(t) \in \partial \Omega$ there exists a unique (different) point $\gamma(t^*)$ such that 
\[
L(\gamma'(t),\gamma'(t^*)) = 0\, .
\]
In other words, given $t\in\mathbb{S}$, there is a unique $t^*\in\mathbb{S}$ ($t^*\neq t$) such that the tangent vectors to $\partial\Omega$ at the points $\gamma(t)$ and $\gamma(t^*)$ are parallel.

\begin{definition}\label{adm}
Given $t_1,t_2\in \mathbb{S}$, we say that $(t_1,t_2)\in \mathbb{S}\times\mathbb{S}$ is \emph{positive admissible} if for a lift $T_1$ of $t_1$ (i.e., $\pi(T_1)=t_1$) and the lift $T_1^*$ of $t_1^*$ such that $T_1<T_1^*<T_1+2\pi $, there exists a lift $T_2$ of $t_2$ (i.e., $\pi(T_2)=t_2$) such that
\[
T_1<T_2<T_1^*.
\]
\end{definition}
\noindent We refer to
\begin{equation} \label{PS}
\hat{\mathcal{P}} = \{ (t_1,t_2) \in \mathbb{S} \times \mathbb{S}: \ (t_1,t_2) \text{ is positive admissible}\}
\end{equation}
as the (open, positive) phase space.

\begin{definition}\label{def sym bill}
The symplectic billiard map of the domain $\partial\Omega$, parametrized by $\gamma\colon \mathbb{S}\to\mathbb{R}^2$, is 
$$\hat{T}: \hat{\mathcal{P}} \rightarrow \hat{\mathcal{P}}, \qquad (t_1,t_2) \mapsto (t_2,t_3)$$
where $\gamma(t_3)\in\partial\Omega$ is the unique point satisfying 
$$L( \gamma^{\prime} (t_2),\gamma(t_3)-\gamma(t_1))=0\, .$$ 
\end{definition}


\noindent The next proposition summarizes the main properties of a symplectic billiard map. We refer to \cite[Section 2]{AT} for further details and for the proofs.

\begin{figure}[ht]
\centering
\includegraphics[scale=0.4]{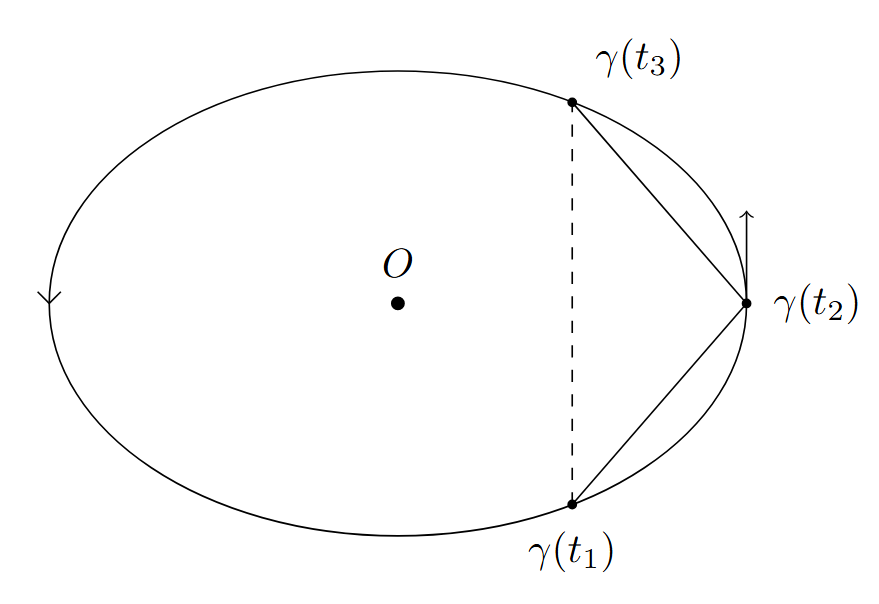}
    \caption{The symplectic billiard map reflection: after the points $\gamma(t_1)$ and $\gamma(t_2)$, the next bounce occurs at the point $\gamma(t_3)$.}
    \label{SB}
\end{figure}

\begin{proposition}\label{properties}
Let $\Omega\subset \mathbb{R}^2$ be a $C^k$ strictly convex smooth domain, $k\geq 2$. Let $\gamma\colon \mathbb{S}\to \mathbb{R}^2$ be a parametrization of $\partial\Omega$. Denote by $\hat T\colon \hat{\mathcal{P}}\to \hat{\mathcal{P}}$ the associated symplectic billiard map. The following properties hold.
\begin{enumerate}
    \item 
    $\hat{T}$ is $C^{k-1}$ and it extends continuously to the closure of $\hat{\mathcal{P}}$ so that 
$$\hat{T}(t,t) = (t,t) \qquad \text{and} \qquad \hat{T}(t,t^*) = (t^*,t)\, .$$
\item For every $(t_1,t_2)\in\hat{\mathcal{P}}$ one has
\begin{equation} \label{GF}
\hat{T}(t_1,t_2) = (t_2,t_3) \qquad \Longleftrightarrow \qquad L_2(t_1,t_2) + L_1(t_2,t_3) = 0
\end{equation}
where we use the notation 
$$L_2(t_1,t_2) := L(\gamma(t_1),\gamma'(t_2)) \qquad \text{and} \qquad L_1(t_2,t_3) := L(\gamma'(t_2),\gamma(t_3))\, .$$


\item The map $\hat T$ does not depend on the choice of the origin $O$. 
\item\label{punto 4 trsnf affini} The map $\hat T$ commutes with any map obtained as affine transformation of the plane, since they preserve tangent directions. 
\end{enumerate}
\end{proposition}


\noindent Up to a change of coordinates, it is possible to see the symplectic billiard map as a (negative) exact twist map (see Definition \ref{aggiunta anna}). In fact the function
\begin{equation}\label{change coordinates}
\phi: \hat{\mathcal{P}} \to \mathbb{S}\times\mathbb{R}, \qquad (t_1,t_2) \mapsto (t_1, -L_1(t_1,t_2))
\end{equation}
is a diffeomorphism onto its image. 
\begin{figure}
    \centering
    \includegraphics[scale=0.15]{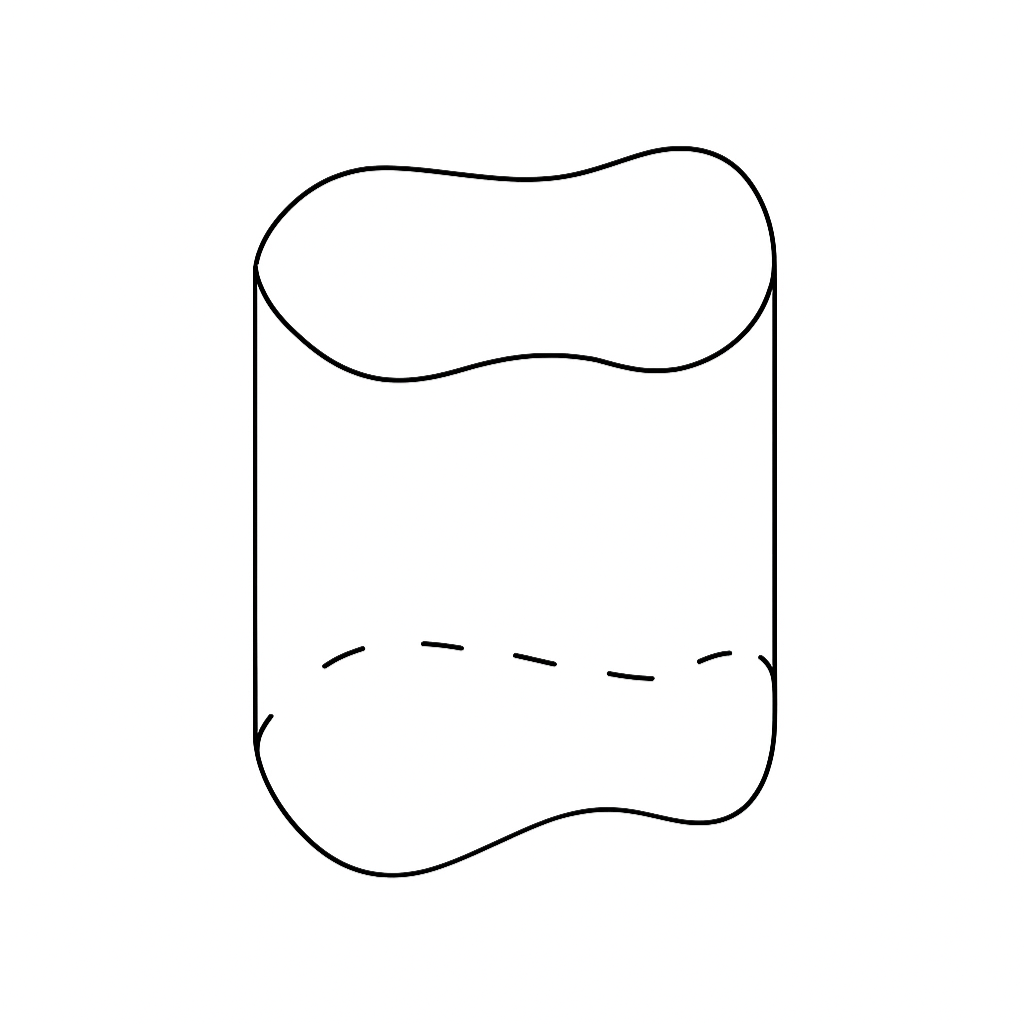}
    \caption{The phase space $\mathcal{P}$.}
    \label{fig:cilindrosbilenco}
\end{figure}The image $\phi(\hat{\mathcal{P}})$ is the open set
\[
\mathcal{P}=\{ (t,s)\in\mathbb{S}\times\mathbb{R} :\ s\in (\psi_1(t),\psi_2(t) )\}\, ,
\]
(see Figure \ref{fig:cilindrosbilenco}), where 
\[\psi_1\colon t\in\mathbb{S}\to \psi_1(t):=-L_1(t,t^*)=-L(\gamma'(t),\gamma(t^*)-O)\in\mathbb{R}\, ,\]
\[
\psi_2\colon t\in\mathbb{S}\to \psi_2(t):=-L_1(t,t)=-L(\gamma'(t),\gamma(t)-O)\in\mathbb{R}\, .
\]
Observe that $\psi_1<0<\psi_2$, so in particular the zero section $\mathbb{S}\times\{0\}$ is contained in $\mathcal{P}$.
By the variational condition \eqref{GF}, denoting by $(t_0,t_1)\in\hat{\mathcal{P}}$ the point such that $\hat T(t_0,t_1)=(t_1,t_2)$, we have also that the second component of $\phi(t_1,t_2)$ equals $L_2(t_0,t_1)$. We can then consider the map
\begin{equation}\label{sympl bill coordinate t s}
T:= \phi\circ \hat T\circ \phi^{-1}\vert_{\mathcal{P}}\colon \mathcal{P}\to\mathcal{P}
\end{equation}
which preserves the area form $dt\wedge ds$. 
Moreover, the map $T$ is a negative twist map, i.e. if $T(t_1,s_1)=(t_2,s_2)$ then
\[
\dfrac{\partial t_2}{\partial s_1}<0\, .
\]
We refer to \cite[Lemma 2.7]{AT} for the proof of this fact. 
\begin{remark}
According to Proposition \ref{prop gen func} and by \eqref{GF}, the function
\[
S\colon (T_1,T_2)\in\mathbb{R}^2\to S(T_1,T_2):= L(\gamma(\pi(T_1)),\gamma(\pi(T_2)))\in\mathbb{R}
\]
is a generating function for the twist map $T$. In the sequel, with an abuse of notation, we will simply write $L(t_1,t_2)$ to refer to $L(\gamma(\pi(T_1)),\gamma(\pi(T_2)))$, where $\pi(T_i)=t_i$ for $i=1,2$; the notations $L_i$ and $L_{ij}$, for $i,j\in\{1,2\}$, will denote the partial derivatives of order one and two, respectively. In the sequel, in order to lighten the notation, we consider also the extension of the function $\gamma$ to $\mathbb{R}$, and so drop $\pi$ in the notation of the generating function.
\end{remark}
\begin{proposition} \label{T dopo T twist}
Let $\Omega\subset \mathbb{R}^2$ be a $C^k$ strictly convex domain, $k\geq 2$. Denote by $T \colon \mathcal{P} \to \mathcal{P}$ the associated symplectic billiard map. The map $T^2$ is a negative twist map.
\end{proposition}
\begin{proof}
Following the same notation as before, we denote $T(t_1,s_1) = (t_2,s_2)$. Then $T^2(t_1,s_1) = (t_3,s_3)$ is a negative twist map if $\frac{\partial t_3}{\partial s_1} < 0$ or, equivalently, 
\begin{equation}\label{bordeaux}
\frac{\partial s_1}{\partial t_3} = -\frac{\partial L_1}{\partial t_3}(t_1,t_2(t_1,t_3))=-L_{12}(t_1,t_2(t_1,t_3))\frac{\partial t_2}{\partial t_3} (t_1,t_3)<0\,.
\end{equation}
Fixed $t_1$, in the above formula, $t_2:=t_2(t_1,t_3)$ gives the unique --by strict convexity-- point such that
\begin{equation} \label{var1}
 L_2(t_1,t_2)+L_1(t_2,t_3)=0\,.
\end{equation}
Condition \eqref{bordeaux} is then easy to verify. In fact, $L_{12}(t_1,t_2)>0$ since $L_{12}(t_1,t_2)=\mathrm{det}(\gamma'(t_1),\gamma'(t_2))$ and $(t_1,t_2)$ is positive admissible, see Definition \ref{adm}. Moreover, considering $t_2=t_2(t_1,t_3)$ and differentiating \eqref{var1} with respect to $t_3$, we obtain
\[
 \big(L_{22}(t_1,t_2)+L_{11}(t_2,t_3)\big)\frac{\partial t_2}{\partial t_3}+L_{12} (t_2,t_3) =0\, .
\]
Observe that $L_{22}(t_1,t_2)+L_{11}(t_2,t_3)=\mathrm{det}(\gamma''(t_2),\gamma(t_3)-\gamma(t_1))$, which is strictly negative since $\gamma''(t_2)$ points into the interior of the domain and $\gamma(t_3)-\gamma(t_1)$ is parallel and co-oriented with $\gamma'(t_2)$. Since also $L_{12}(t_2,t_3)$ is positive, we deduce that
\[
\frac{\partial t_2}{\partial t_3}=-\frac{L_{12} (t_2,t_3)}{L_{11}(t_2,t_3)+L_{22}(t_1,t_2)}>0\, .
\]
Thus, \eqref{bordeaux} is true. This concludes the proof.
\end{proof}


\subsection{Dissipative symplectic billiards}\label{INIZIO}

Let $T: \mathcal{P} \ni (t_1,s_1) \mapsto (t_2,s_2) \in \mathcal{P}$ be the symplectic billiard map on a strictly convex domain $\Omega$ with $C^k$ boundary, $k\geq 2$. Let $\gamma\colon\mathbb{S}\to\mathbb{R}$ be a parametrization of the boundary $\partial\Omega$ and let $O\in \mathrm{int}(\Omega)$ be the fixed origin. Let $\lambda \in (0,1)$: it will be our fixed dissipative parameter. Denote by
\[\mathcal{H}_{\lambda}: \mathcal{P} \to \mathcal{P}, \qquad \mathcal{H}_{\lambda}(t_1,s_1) = (t_1,\lambda s_1)\]
the $\lambda$-contraction map along the fiber. Observe that $\mathcal{H}_\lambda(\mathcal{P})\subset \mathrm{int}(\mathcal{P})$.
\begin{definition} \label{disss}
The dissipative symplectic billiard map on $\Omega$ is defined as 
$$T_{\lambda} := \mathcal{H}_{\lambda} \circ T: \mathcal{P} \to \mathcal{P}, \qquad T_{\lambda}(t_1,s_1) = (t_2,\lambda s_2)\, .$$
\end{definition}
\noindent Since $T$ preserves the standard area form $dt\wedge ds$ and $T$ is a twist map, it can be shown that $T_{\lambda}$ is a dissipative, twist map. In particular, $T_\lambda$ dissipates the area of a constant factor $\lambda$:
\[T^*_{\lambda} (dt \wedge ds) = T^*(\mathcal{H}_{\lambda}^*(dt \wedge ds)) = \lambda dt \wedge ds\, .\]

According to the literature (see \cite{AHV,AA,AF,AFR,MS}), we can also say that the map $T_\lambda$ is a conformally symplectic map: the symplectic form is contracted by the dynamics by some factor. In dimension $2$, symplectic forms are area forms.
\noindent We will be interested in studying the (Birkhoff) attractor of such dissipative, symplectic billiard maps. A straightforward consequence of the area dissipation is the next lemma.
\begin{lemma}\label{Luca}
Let $\lambda\in(0,1)$. Consider the dissipative symplectic billiard map $T_\lambda$, introduced in Definition \ref{disss}. Denoting by $\Lambda$ its Birkhoff attractor, we have
    \[\Lambda\cap\left(\mathbb{S}\times\{0\}\right)\neq\emptyset\, .\]
\end{lemma}
\begin{proof}
    In the proof, we indicate $\mathcal{P}^+ := \{(t,s)\in\mathcal{P}: s>0\}$ and $\mathcal{P}^- := \{(t,s)\in\mathcal{P}: s<0\}$. Suppose by contradiction that $\Lambda\cap\left(\mathbb{S}\times\{0\}\right)=\emptyset$. Then, without loss of generality and since $\mathbb{S}\times\{0\}\subset \mathcal{P}$, we can assume that \begin{equation}\label{birk tutto sopra}\Lambda \subset\mathcal{P}^+\, .
    \end{equation}
    The Birkhoff attractor separates the annulus: denote then $\mathcal{P}\setminus\Lambda=U^\Lambda_-\sqcup U^\Lambda_+$. 
    Since $\Lambda$ is $T_\lambda$-invariant and by the definition of $T_\lambda$, we deduce that $T(\Lambda)=\mathcal{H}_{\frac{1}{\lambda}}\circ T_\lambda(\Lambda) = \frac{1}{\lambda} \Lambda\subset \mathcal{P}^+$.
    From this observation, by \eqref{birk tutto sopra} and since $T$ preserves boundaries, we deduce that  
$$m(T(U^\Lambda_-)) = m(\mathcal{P}^-) + \frac{1}{\lambda} m(U^\Lambda_- \setminus \mathcal{P^-}) > m(U^-_\Lambda)\, ,$$ which contradicts the conservative properties of $T$.
\end{proof}
\noindent According to \eqref{change coordinates}, let 
\begin{equation} \label{cambio}
\phi^{-1}\colon  \mathcal{P} \to \hat{\mathcal{P}}, \qquad \phi^{-1}(t_1, s_1) = \phi(t_1,-L_1(t_1,t_2)) = (t_1,t_2)\, .
\end{equation}
Observe that it remains well-defined on the image of $T_\lambda$, since the map is dissipative. By $\mathcal{O}(t_1,s_1)$ and $\mathcal{O}(t_1,t_2)$ we denote the orbit of $(t_1,s_1)$ under $T_\lambda$ and of $(t_1,t_2)$ under $\phi^{-1}\circ T_\lambda\circ \phi$, respectively.

The following lemma provides a formula for the differential of $T_\lambda$, which will be useful later on.

\begin{lemma}\label{formula DT lambda}
    Let $T_\lambda$ be a dissipative symplectic billiard map. Let $(t_1,s_1)\in\mathcal{P}$ and let $(t_1,t_2)=\phi^{-1}(t_1,s_1)\in\hat{\mathcal{P}}$. Then
    \begin{equation}
\label{differential}
DT_\lambda (t_1,s_1) = -\frac{1}{L_{12}(t_1,t_2)}\begin{pmatrix}
      L_{11}(t_1,t_2) & 1\\
-\lambda L_{12}^2(t_1,t_2)+\lambda L_{22}(t_1,t_2)\cdot L_{11}(t_1,t_2)& \lambda L_{22}(t_1,t_2)
    \end{pmatrix}\, ;
\end{equation}
in particular, we have that $\det DT_\lambda (t_1,s_1) = \lambda$.
\end{lemma}
\begin{proof}
    
\noindent From Definition \ref{disss}, we have that $T_\lambda=\mathcal{H}_\lambda\circ T$; thus, for every $(t_1,s_1)\in\mathcal{P}$,
\begin{equation}\label{formula diff 1}
DT_\lambda(t_1,s_1)=D\mathcal{H}_\lambda(T(t_1,s_1))DT(t_1,s_1)=\begin{pmatrix}
    1 & 0 \\ 0 & \lambda
\end{pmatrix}\, DT(t_1,s_1)\, .
\end{equation}
Since $T$ is conservative, it follows immediately that $\det DT_\lambda(t_1,s_1)=\lambda$. Considering $(t_1,t_2)=\phi^{-1}(t_1,s_1)$, we have that
\[
T\circ \phi (t_1,t_2) = (t_2, L_2(t_1,t_2))=(t_2,L(\gamma(t_1),\gamma'(t_2))\, ,
\]
and so
\[
D(T\circ \phi)(t_1,t_2) =\begin{pmatrix}
    0 & 1 \\ \frac{\partial}{\partial t_1}L_2(t_1,t_2) & \frac{\partial}{\partial t_2}L_2(t_1,t_2)
\end{pmatrix}=\begin{pmatrix}
    0 & 1 \\ L_{12}(t_1,t_2) & L_{22}(t_1,t_2)
\end{pmatrix}\, .
\]
Since $DT(t_1,s_1)=D(T\circ\phi)(t_1,t_2)\left(D\phi(t_1,t_2)\right)^{-1}$, and $\left(D\phi(t_1,t_2)\right)^{-1}=\frac{1}{-L_{12}(t_1,t_2)}\begin{pmatrix}
    -L_{12}(t_1,t_2) & 0 \\ L_{11}(t_1,t_2) & 1
\end{pmatrix}$, we conclude that
\begin{equation}\label{formula diff 2}
DT(t_1,s_1)= \begin{pmatrix}
    0 & 1 \\ L_{12}(t_1,t_2) & L_{22}(t_1,t_2)
\end{pmatrix}\, \begin{pmatrix}
    1 & 0 \\ -\frac{L_{11}(t_1,t_2)}{L_{12}(t_1,t_2)} & -\frac{1}{L_{12}(t_1,t_2)}\, .
\end{pmatrix}
\end{equation}
From \eqref{formula diff 1} and \eqref{formula diff 2}, we conclude that
\[
DT_\lambda(t_1,s_1)= -\frac{1}{L_{12}(t_1,t_2)}\begin{pmatrix}
    L_{11}(t_1,t_2) & 1 \\ -\lambda L_{12}^2(t_1,t_2) +\lambda L_{11}(t_1,t_2)L_{22}(t_1,t_2) & \lambda L_{22}(t_1,t_2)
\end{pmatrix}\, ,
\]
as stated.
\end{proof}


\noindent In the next lemma, we proceed by giving a geometrical characterization of Definition \ref{disss}. According to Definition \ref{def sym bill} and to \eqref{sympl bill coordinate t s}, we denote $\hat{T}_\lambda \colon \hat{\mathcal{P}}\to\hat{\mathcal{P}}$ the map $\phi^{-1}\circ T_\lambda\circ\phi$.
\begin{lemma} \label{origine} Let $\gamma: \mathbb{S} \ni t \mapsto \gamma(t) \in \partial \Omega$ be a parametrization of $\partial \Omega$. Then $\hat{T}_{\lambda}(t_1,t_2) = (t_2,t_3)$ if and only if $\gamma(t_3) - \lambda \gamma(t_1) \in T_{\gamma(t_2)} \partial \Omega$.
\end{lemma}
\begin{proof}
Fixed $(t_1,t_2) \in \mathbb{S} \times \mathbb{S}$, we indicate $\hat{T}(t_1,t_2) = (t_2,\tilde{t}_3)$ and $\hat{T}_{\lambda}(t_1,t_2) = (t_2,t_3)$. Equivalently, in the $(t,s)$ coordinates:
$$T(t_1,s_1) = (t_2,s_2) \qquad \text{with } s_2 = -L_1(t_2,\tilde{t}_3)$$
and
$$T_{\lambda}(t_1,s_1) = (t_2,\lambda s_2) \qquad \text{with } \lambda s_2 = -L_1(t_2,t_3)\, .$$
Together with the variational condition \eqref{GF} for $\hat{T}$, that is $L_2(t_1,t_2) + L_1(t_2,\tilde{t}_3) = 0$, we obtain that
\[\lambda L_2(t_1,t_2) + L_1(t_2,t_3) = 0 \Leftrightarrow L(\lambda\gamma(t_1), \gamma'(t_2)) + L(\gamma'(t_2),\gamma(t_3)) = L(\gamma'(t_2), -\lambda\gamma(t_1)+\gamma(t_3))=0\, ,\]
and the statement follows. 
\end{proof}

\begin{remark}
    The definition of $T_\lambda$ depends on the choice of the origin $O$, which is not the case for the (conservative) symplectic billiard map $T$.
\end{remark}
In the sequel, we discuss the role of $4$-periodic points in order to motivate the choice of the origin we will made in the next sections. 

\begin{lemma}\label{origin compatible}
    Let $T$ be a symplectic billiard map. Then, there exists a choice of the origin $O$ such that the set $\Pi$ of 4-periodic points whose orbit is contained in the zero section is not empty. We will say that the origin $O$ is \emph{compatible}.
\end{lemma}

\begin{proof}
Being $T$ a conservative twist map,
by a classical Birkhoff's theorem \cite{B22}, $T$ possesses at least two $4$-periodic orbits. Let $\{ t_i\}_{i = 1}^4$ be the set of points in $\mathbb{S}$ corresponding to one of such orbits. 
In particular, by the definition of $\hat T$, we have
\[
L(\gamma'(t_2),\gamma(t_3)-\gamma(t_1))=L(\gamma'(t_4),\gamma(t_1)-\gamma(t_3))=0
\]
and
\[
L(\gamma'(t_3),\gamma(t_4)-\gamma(t_2))=L(\gamma'(t_1),\gamma(t_2)-\gamma(t_4))=0\, ;
\]
this implies then
\[
L(\gamma'(t_1),\gamma'(t_3))=L(\gamma'(t_2),\gamma'(t_4))=0\, ,
\]
that is, the vectors $\gamma'(t_1)$ and $\gamma'(t_3)$ are parallel, as well as the vectors $\gamma'(t_2)$ and $\gamma'(t_4)$.
Consequently, in the phase space $\hat{\mathcal{P}}$, a
$4$-periodic orbit for $\hat T$ is given by
$$
\{ (t_1,t_2), (t_2,t_1^*), (t_1^*, t_2^*), (t_2^*, t_1) \}\, .  
$$
Consider then the (inscribed) quadrilateral with vertices in $\{ \gamma(t_i)\}_{i = 1}^4$. Fix the origin $O\in \mathrm{int}(\Omega)$ as the intersection of the diagonals of the inscribed quadrilateral corresponding to (at least) one of the $4$-periodic orbits for $\hat T$. As a consequence, by this choice, in the phase space $\mathcal{P}$, any such $4$-periodic orbit is given by
\begin{equation} \label{quattro orbita}
\{ (t_1,0), (t_2,0), (t_1^*,0), (t_2^*,0)\}\, ,
\end{equation}
i.e., its $T$-orbit is contained in the zero section, as required. This follows by the definition of the change of coordinates $\phi$ in \eqref{change coordinates}, since $\phi(t_1,t_2)=(t_1, -L(\gamma'(t_1),\gamma(t_2)-O))$, and from the fact that $L(\gamma'(t_1),\gamma(t_2)-O)=0$ since $\gamma(t_2)-O$ belongs to the line generated by $\gamma(t_2)-\gamma(t_4)$, by the choice of the origin, and since $\gamma(t_2)-\gamma(t_4)$ is parallel to $\gamma'(t_1)$.
\end{proof}

An immediate consequence is the following one.
\begin{corollary}\label{O compatible then points survive}
Let $\Omega$ be a $C^k$ strictly convex domain, $k\geq 2$. Choose a compatible origin $O$. Let $T$ be the conservative symplectic billiard map and, for every $\lambda\in(0,1)$, let $T_\lambda$ be the associated dissipative symplectic billiard map.  For every $(t,0)\in\Pi$, one has
\begin{equation}\label{no effect}
T_\lambda(t,0)=T(t,0)\, .
\end{equation}
\end{corollary}


We will see in the following remark that, for a centrally symmetric domain, the natural choice of the origin as the center of symmetry is compatible.

\begin{remark} \label{tutte} Let $\Omega$ be a (strictly convex) centrally symmetric domain, i.e., there is a point $O\in\mathbb{R}^2$ such that $\partial\Omega$ is invariant under the isometric involution $(x-O,y-O)\mapsto (-x-O,-y-O)$. The point $O$ is called the center of symmetry. Fix it as the origin.
Then, the phase space $\mathcal{P}$ becomes
\[\mathcal{P} = \{ (t,s) \in \mathbb{S} \times \mathbb{R}: \ -\psi(t) < s < \psi(t)\}\, ,\]
where $\psi(t) = -L_1(t,t)=-L(\gamma'(t),\gamma(t)-O)$. In such a case, for every $t\in\mathbb{S}$, the corresponding $t^*\in\mathbb{S}$ is just $t+\pi$. Thus, any $4$-periodic orbit has the form
\begin{equation*} \label{serve dopo}
\{ (t_1,t_2), (t_2,t_1 + \pi), (t_1 + \pi, t_2 + \pi), (t_2 + \pi,t_1) \}\, .
\end{equation*}
By this choice, it can be proved that the orbit of any $4$-periodic point is contained in the zero section $\mathbb{S}\times\{0\}$, i.e., the set $\Pi$ is the set of all $4$-periodic points.
\end{remark}

As a concluding example, we show now that the Birkhoff attractor for dissipative symplectic billiard maps on a (centrally symmetric) Radon domain is the zero section 
$\mathbb{S} \times \{ 0\}$. A Radon domain $\Omega$ is a centrally symmetric domain such that every point of the boudary $\partial\Omega$ is the vertex of an inscribed parallelogram of maximal area, see \cite[Section 3]{MaSw}. 


\begin{proposition} \label{Radon doamins}
Let $\Omega$ be a Radon domain and let the origin $O$ be the center of symmetry. Then, for every $\lambda\in(0,1)$, the Birkhoff attractor of $T_\lambda$ is $\Lambda=\mathbb{S}\times\{0\}$.
\end{proposition}
\begin{proof} By the property of centrally symmetric Radon domains recalled above, for the corresponding symplectic map $T$, the zero section $\mathbb{S}\times\{0\}$ is an invariant curve made up of $4$-periodic points. By Corollary \ref{O compatible then points survive}, since he chosen origin is compatible, for every $\lambda\in(0,1)$, the zero section is also $T_\lambda$-invariant. By minimality with respect to the inclusion (see Proposition \ref{Birk-attr-minimal}), we conclude that $\mathbb{S}\times\{0\}$ is the Birkhoff attractor for $T_\lambda$. 
\end{proof}

\begin{remark} \label{RADON and ELLIPSE}
Let $\Omega$ be an elliptic domain. In particular, it is Radon: by Proposition \ref{Radon doamins}, for every $\lambda$, the Birkhoff attractor is the zero section. Moreover, it always coincides with the attractor $\Lambda_0$, as we will argue now.
As proved in \cite{BaBe}, elliptic domains are the only completely integrable ones with respect to the symplectic reflection law. Up to an affine transformation, we can assume that $\Omega$ is the unitary circle.
Considering as first coordinate $\alpha\in\mathbb{S}$ the counterclockwise oriented angle with the fixed direction $(1,0)$, the symplectic billiard map becomes
$
T(\alpha,s)=(\alpha+\arccos{(s)},s)$. Thus, the second coordinate is a first integral.
For every $\lambda\in(0,1)$, the dissipative symplectic billiard map is $T_\lambda(\alpha,s)=(\alpha+\arccos{s},\lambda s)$; in particular, the second coordinate is now a Lyapunov function. The corresponding neutral set, which is the attractor $\Lambda_0$, coincides with the zero section, i.e., with the Birkhoff attractor $\Lambda$. It could be interesting to construct a Radon domain (hence different to an ellipse) for which $\Lambda_{\lambda} \subsetneq \Lambda_{\lambda}^0$.

\end{remark}

\begin{remark} 
For dissipative Birkhoff billiards, elliptic domains play a special role: it is possible to completely described their Birkhoff attractor for any dissipative factor $\lambda\in(0,1)$ (see \cite[Theorem 4.6]{BFL}). Given an elliptic table with non-zero eccentricity, the $2$-periodic orbits --corresponding to the trajectories along the major and minor axes-- are hyperbolic orbits of saddle and sink type respectively, and they are not affected by dissipation. The Birkhoff attractor $\Lambda$ coincides with the global one $\Lambda_0$, and it is the closure of the unstable manifold for the $2$-periodic hyperbolic orbit. 

When the dissipation is strong enough and the eccentricity is small enough, the Birkhoff attractor turns out to be a normally contracted graph over $\mathbb{S} \times \{0\}$, see \cite[Theorem 4.6, Proposition 5.16]{BFL}. It is worth noting that the proof of \cite[Theorem 4.6]{BFL} relies on the existence of a first integral for the conservative elliptic Birkhoff billiard dynamics. In fact, such an integral becomes a Lyapunov function when the dissipation is turned on. Its neutral set is the set of $2$-periodic points. The global attractor, which corresponds to the Birkhoff one, is then detected by using such a Lyapunov function.

For dissipative symplectic billiards, it is not known if there exists a class of tables which can replace the class of elliptic domains for dissipative Birkhoff billiards. The main reason is that it is not known if there exists a class of tables, for (conservative) symplectic billiards, whose dynamics is integrable, but not totally integrable.

\end{remark}

\color{black}

\section{Normally contracted Birkhoff attractors} \label{NORMALLY}
In Proposition \ref{Radon doamins},  
we have seen that the Birkhoff attractor of any Radon billiard table is very simple, no matter the strength of the dissipation: it is always $\Lambda = \mathbb{S}\times \{0\}$. The aim of this section is showing that --when the dissipation is strong, i.e., $0<\lambda\ll 1$ -- the corresponding dissipative symplectic billiard map of a strongly convex $C^k$ domain, $k\geq 2$, 
exhibits a topologically simple Birkhoff attractor, i.e., $\Lambda$ is a normally contracted graph and it coincides with the attractor $\Lambda_0$. This section follows very closely \cite[Section 5]{BFL}. \\
~\newline
The next result is the main point to obtain a cone-field criterion and then invoke standard results in normally hyperbolic dynamics, recalled later on. 
In the next proposition, as in Theorem \ref{nc}, the choice of the origin $O \in \mathrm{int}(\Omega)$ is arbitrary. 
\begin{proposition}\label{cones}
 Let $\Omega$ be a $C^k$ strongly convex domain, $k\geq 2$. Then, there exist $\lambda(\Omega) \in(0,1)$, $M > 0$, $\alpha>0$ and a cone-field $\mathcal{C}_{\alpha}=(\mathcal{C}_{\alpha}(t,s))_{(t,s)\in\mathcal{P}}$ containing the horizontal direction
\[\mathcal{C}_{\alpha}(t,s)=\lbrace v\in T_{(t,s)}\mathcal{P}: v=(v_t,v_s), |v_s|\leq\alpha|v_t|\rbrace\, ,\]
such that, for every $\lambda\in(0,\lambda(\Omega))$ and for every $(t,s)\in\mathbb{S}\times[-M \cdot \lambda(\Omega),M \cdot \lambda(\Omega)]$, one has
\[DT_\lambda(t,s)\mathcal{C}_{\alpha}(t,s)\subset \mathrm{int}(\mathcal{C}_{\alpha}(T_\lambda(t,s)))\cup\lbrace0\rbrace\, .\]

\end{proposition}
\begin{proof}
Consider the arc-length parametrization $\gamma\colon\mathbb{S}\to \partial\Omega$ of the billiard table $\Omega$. Define the following positive constants
\[
C_1:=\max_{t_1,t_2\in \mathbb{S}}\vert L_{11}(t_1,t_2)\cdot L_{22}(t_1,t_2)-L^2_{12}(t_1,t_2)\vert\, ,
\]
\[C_2:=\max_{t_1,t_2\in\mathbb{S}}\vert L_{22}(t_1,t_2)\vert
\]
and
\[M(\Omega) :=\max\biggl\lbrace\max_{t\in\mathbb{S}}|L_1(t,t)|,\max_{t\in\mathbb{S}}|L_1(t,t^*)| \biggr\rbrace\, .\]
Since $\gamma$ is the arc-length parametrization, for every $t\in\mathbb{S}$, one has that $\gamma''(t)=k(t)\, i\gamma'(t)$, where $k(t)\in\mathbb{R}$ denotes the curvature of the table at the point $\gamma(t)$, and the multiplication by $i$ is the rotation of angle $\frac{\pi}{2}$, so that $\gamma''(t)$ always points towards the inside of the domain. Recall that $L(\gamma(t_1),\gamma(t_2))$ refers to the determinant of the matrix whose columns are, respectively, $\gamma(t_1)-O$ and $\gamma(t_2)-O$.

By continuity of the involved function and compactness, we can choose $\lambda_1\in(0,1)$ such that there exists $c_0>0$ so that for every $(t_1,s_1)\in\mathcal{P}$, if
\[|s_1|=|L_1(t_1,t_2)|=\vert \mathrm{det}(\gamma'(t_1),\gamma(t_2)) \vert<M(\Omega) \cdot \lambda_1 \, ,\]
that is, if $\gamma'(t_1)$ is not so far from being parallel to $\gamma(t_2)-O$, then
\[\vert L_{11}(t_1,t_2)\vert=\vert 
L(\gamma''(t_1),\gamma(t_2) \vert=\vert k(t_1)\vert \, \vert L(i\gamma'(t_1),\gamma(t_2))\vert \geq c_0>0\, ,\]
since $i\gamma'(t_1)$ is then not so far from being perpendicular to $\gamma(t_2)-O$.
Fix some $0<\alpha<c_0$. 
Let $\mathcal{C}_{\alpha}$ be the cone-field defined as \[\mathcal{C}_{\alpha}(t,s)=\lbrace v\in T_{(t,s)}\mathcal{P}: v=(v_t,v_s), |v_s|\leq\alpha|v_t|\rbrace\, .\]
We are considering on each tangent space the coordinates inherited from $(t,s)\in\mathcal{P}$. 
By Lemma \ref{formula DT lambda}, for every  $(t,s)\in 
\mathbb{S}\times[-M(\Omega) \cdot \lambda_1,M(\Omega) \cdot \lambda_1]$ and for every vector $v=(a,b)\in\mathcal{C}_{\alpha}(t,s)$, we get 
  \begin{equation*}
  \begin{split}
  v':=&DT_\lambda(t,s)v\\=&-\frac{1}{L_{12}(t_1,t_2)}\begin{pmatrix}
       aL_{11}(t_1,t_2)+b\\
        \lambda a\left[L_{11}(t_1,t_2)\cdot L_{22}(t_1,t_2)-L^2_{12}(t_1,t_2)\right]+\lambda b L_{22}(t_1,t_2)
    \end{pmatrix}\\:=&-\frac{1}{L_{12}(t_1,t_2)}\begin{pmatrix}
        a'\\
        b'
    \end{pmatrix}.
    \end{split}
  \end{equation*}
  Thus we have
  \begin{equation*}
  \begin{split}
      |a'|&=| aL_{11}(t_1,t_2)+b|\geq|a||L_{11}(t_1,t_2)|-|b|\geq|a|(c_0-\alpha),\\
      |b'|&=| \lambda a\left[L_{11}(t_1,t_2)\cdot L_{22}(t_1,t_2)-L^2_{12}(t_1,t_2)\right]+\lambda b L_{22}(t_1,t_2)|\\
      &\leq\lambda |a|(C_1+\alpha C_2)\, .
        \end{split}
  \end{equation*}
  Given now $\mu_0\in (0,1)$, it holds
 \[DT_\lambda(t,s)\mathcal{C}_{\alpha}(t,s) \subset \mathcal{C}_{\mu_0\alpha}(T_\lambda(t,s))\subset \mathrm{int}(\mathcal{C}_{\alpha}(T_\lambda(t,s)))\cup\{0\}\, ,\]
 where $\mathcal{C}_{\mu_0\alpha}(t,s):=\{v\in T_{(t,s)}\mathcal{P} :\ v=(v_t,v_s), \vert v_s\vert\leq \mu_0\alpha\vert v_t\vert\}$, provided that $\lambda\in(0,\lambda(\Omega))$, with
 \begin{equation*}
     \lambda(\Omega):=\text{min}\left\{\lambda_1,\frac{\mu_0\alpha(c_0-\alpha)}{K+\alpha C_0}\right\}\, .
 \end{equation*}
\end{proof}

The rest of the section follows \cite[Section 5]{BFL}: we are able to show that, when the dissipation is strong, the Birkhoff attractor of the corresponding dissipative symplectic map coincides with the attractor and it is a \emph{normally contracted} graph. Even if the ideas are the same as in \cite{BFL}, for sake of clarity, we recall here the main definitions, coming from hyperbolic dynamics, we state the main results and give an idea of the proof.

\begin{remark}
    As we will see in the main statements, the results of this section hold for \textbf{every} $C^k$ strongly convex domain, with $k\geq 2$, while, in the Birkhoff dissipative dynamics, the analogous results hold for a class of strictly convex domains satisfying some geometric \textit{pinching} condition (see \cite[Definition D]{BFL}).
\end{remark}

Let us start by recalling some well-known definitions and results in normally hyperbolic dynamics. We refer to \cite{BergerBounemoura,CroPot, HPS, Sambarino}.
\begin{definition}
    Let $M$ be a compact Riemmannian manifold without boundary and $f\colon M\to M$ be a $C^k$ diffeomorphism, with $k\geq 1$. A compact invariant set $K$ has a dominated splitting if $T_KM=E\oplus F$, where the $Df$-invariant continuous subbundles $E$ and $F$ are such that there exists $C>0$ and $\nu\in (0,1)$ so that, for every $x\in K$
    \[
    \Vert Df^n(x)\vert_E\Vert \cdot \Vert Df^{-n}(f^n(x))\vert_F\Vert \leq C\nu^n \qquad \forall n\geq 0\, .
    \]
\end{definition}
\begin{definition}
    Let $M$ be a compact Riemmannian manifold without boundary and $f\colon M \to M$ be a $C^k$ diffeomorphism, with $k\geq 1$. Let $N$ be a closed $C^1$ $f$-invariant manifold. Then $N$ is normally contracted if $N$ has a dominated splitting $T_NM=E^s\oplus TN$ such that $E^s$ is uniformly contracted, i.e., there exists $n_0\in\mathbb{N}$ and $\mu\in(0,1)$ such that for any $n\geq n_0$ one has $\Vert Df^n(x)\vert_{E^s}\Vert\le \mu^n$ for every $x\in N$.

    We say that $N$ is $k$-normally contracted if $N$ is normally contracted and there exists $C>0$ and $\nu\in(0,1)$ such that for every $x\in N$ and for every $0\le j\le k$ one has $\Vert Df^n(x)\vert_{E^s}\Vert \cdot \Vert Df^{-n}(f^n(x))\vert_{TN}\Vert^j\le C \nu^n$ for all $n\ge 0$.
\end{definition}

Proposition \ref{cones} represents the main point to obtain a cone-field criterion and deduce the following corollary.

\begin{corollary}\label{lambda'}
Let $\Omega$ be a strongly convex domain with $C^k$ boundary, $k \ge 2$. Let $\lambda(\Omega)\in(0,1)$ be given by Proposition \ref{cones}. Then, for every $\lambda\in(0,\lambda(\Omega))$, the attractor $\Lambda_0$ has a dominated splitting $E^s\oplus E^c$ with $E^s$ uniformly contracted. Moreover, every $(t,s)\in\Lambda_0$ has a stable manifold $\mathcal{W}^s(t,s)$ that is transversal to the horizontal. Furthermore, there exists $0<\lambda'(\Omega)<\lambda(\Omega)$ such that for some constant, $C>0$ and $0<\nu<1$, for every $\lambda\in(0,\lambda'(\Omega))$, for every $x\in\Lambda_0$ and for every $1\leq j\leq k-1$
\begin{equation}
    \| DT_\lambda^n(x)\vert_{E^s}\|\cdot\|DT_\lambda^{-n}(T_\lambda^n(x))\vert_{E^c}\|^j\leq C\nu^n,\quad \forall n\geq0\, .
\end{equation}
\end{corollary}
\begin{proof}
The proof follows from the application of the cone-field criterion (see \cite[Theorem 2.6]{CroPot}), and it is \textit{verbatim} the proof of \cite[Proposition 5.5]{BFL}.
\end{proof}
\begin{remark}
    If it would be possible to say \emph{a priori} that $\Lambda_0$ is a $C^1$ manifold, then the previous corollary would say that $\Lambda_0$ is $l$-normally contracted.
\end{remark}
\begin{theorem}\label{nc}
 Let $\Omega$ be a strongly convex domain with $C^k$ boundary, $k \ge 2$. 
 \begin{itemize}
\item[$(a)$] Let $\lambda(\Omega)\in(0,1)$ given in Proposition \ref{cones}. Then for $\lambda\in(0,\lambda(\Omega))$, the Birkhoff attractor $\Lambda$ coincides with the attractor $\Lambda_0$ and it is a normally contracted $C^1$ graph on $\mathbb{S}\times\lbrace0\rbrace$. 
\item[ $(b)$] Let $\lambda'(\Omega)<\lambda(\Omega)$ given in Corollary \ref{lambda'}. Then for $\lambda\in (0, \lambda'(\Omega))$, $\Lambda=\Lambda_0$ is a $C^{k-1}$ graph and $\Lambda$ converges to $\mathbb{S}\times\lbrace0\rbrace$, as $\lambda\to0$, in the $C^1$ topology.
 \end{itemize}
\end{theorem}
We give here the main ideas of the proof and refer to the proof of \cite[Theorem 5.7]{BFL} for full details.
\begin{proof}[Idea of the proof.]The proof is \textit{verbatim} the proof of \cite[Theorem 5.7]{BFL}. Let us give here the main ideas. Consider $\lambda\in(0, \lambda(\Omega))$ and $(\mathcal{C}_{\alpha}(t,s))_{(t,s)\in\mathcal{P}}$ given by Proposition \ref{cones}; define the set of graphs $$\mathcal{F}:=\lbrace\gamma:\mathbb{S}\to\bigl[-\lambda,
\lambda
\bigr] \text{ s.t. }\gamma\in C^1(\mathbb{S})\text{ and }(1,\gamma'(t))\in\mathcal{C}_{\alpha}(t,s),\forall t\in\mathbb{S}\rbrace.$$
Recall that $p_1,p_2$ are the projections on the first and the second coordinate, respectively. By mean of the billiard map $T_\lambda$, it is possible to construct the graph transform
\begin{equation*}
\mathcal{G_{T_\lambda}}:\,\mathcal{F}\to\mathcal{F},\qquad 
        \gamma\mapsto(t,p_2\circ T_{\lambda}(g_\lambda^{-1}(t),\gamma(g_\lambda^{-1}(t))))\, ,
\end{equation*}
where \begin{equation*}
g_\lambda:\,\mathbb{S}\to\mathbb{S},\qquad t\mapsto p_1\circ T_\lambda(t,\gamma(t))\, .
\end{equation*}
In order to obtain $\Lambda=\Lambda_0$, it is sufficient to prove that $\mathcal{G}_{T_\lambda}$ is a contraction on $\mathcal{F}$ for the norm $\|\cdot \|_\infty$. Its only fixed point corresponds precisely to the global attractor and, since it is a graph, it coincides with the Birkhoff attractor. The regularity is a direct consequence of the previous corollary, and the convergence comes from the constructed cone-field.
\end{proof}

\section{Centrally symmetric dissipative billiards} \label{CSSB}

In this section, we consider dissipative symplectic billiard maps for a centrally symmetric domain $\Omega$, i.e., a domain $\Omega$ such that there exists $O\in\mathrm{int}(\Omega)$ for which $\Omega$ is invariant by involution $(x-O,y-O)\mapsto (-x-O,-y-O)$. Under the assumption of strong dissipation, not only Theorem \ref{nc} holds, i.e., the Birkhoff attractor is a normally contracted graph, but we will prove the following results.
\begin{itemize}
\item[$(a)$] for $\lambda \in (0,1)$ small enough, the Birkhoff attractor 
intersects the zero section exactly in the 4-periodic points, that is $(t_1,0) \in \Pi$ if and only if $(t_1,0) \in \Lambda_{\lambda}$; 
\item[$(b)$] for $\lambda \in (0,1)$ small enough, there exists an open and dense set of centrally symmetric domains, whose associated Birkhoff attractor 
has rotation number $1/4$ and decomposes as 
\[\Lambda=\bigcup_{i=1}^l \bigcup_{j=0}^3\overline{\mathcal{W}^u(T_\lambda^j(H_i);T_\lambda^4)}\, ,\]
 where $\lbrace H_i\rbrace_{i=1}^l$ is a finite family of $4$-periodic points of saddle type, and $\mathcal{W}^u(H;T_\lambda^4)$ is the unstable manifold of a saddle point $H$, with respect to the dynamics of $T_\lambda^4$.
 \end{itemize}
 
Recall that $\hat{T}\colon \hat{\mathcal{P}}\to\hat{\mathcal{P}}$ denotes the symplectic billiard map in the coordinates $(t_1,t_2)\in\hat{\mathcal{P}}$. 
Let us introduce the maps 
\begin{align}\label{i}
\hat{I}: \hat{\mathcal{P}} \to \hat{\mathcal{P}}, \qquad \hat{I} (t_1,t_2) = (t_2^*,t_1)\, , \\
\hat{I}_2: \hat{\mathcal{P}} \to \hat{\mathcal{P}}, \qquad \hat{I}_2 (t_1,t_2) = (t_1,t_0^*)\, .
\end{align}
Notice also that 
$$\hat{I}^2(t_1,t_2) = \hat{I} \circ \hat{I} (t_1,t_2) = (t_1^*,t_2^*).$$ 
The corresponding maps $I$ and $I_2$ in $\mathcal{P}$ can be described respectively as
\[I\circ \phi(t_1,t_2)=I(t_1, -L_1(t_1,t_2))=(t_2^*,-L_1(t_2^*,t_1))\]
and
\[I_2\circ\phi(t_1,t_2)=I_2(t_1,-L_1(t_1,t_2))=(t_1,-L_1(t_1,t_0^*))\, .\]
\begin{lemma} \label{secondo}
Let $\Omega$ be a strictly convex $C^k$ domain, $k\ge 2$. Let $O \in \mathrm{int}(\Omega)$ be compatible (see Lemma \ref{origin compatible}). The following properties hold.
\begin{enumerate}
\item\label{punto 1} Recall that $\mathcal{P}=\{(t,s)\in\mathbb{S}\times\mathbb{R}:\ \psi_1(t)\leq s\leq \psi_2(t)\}$. Then, $I_2(\mathrm{graph}(\psi_1))=\mathrm{graph}(\psi_2)$ and $I_2(\mathrm{graph}(\psi_2))=\mathrm{graph}(\psi_1)$. In particular, if $\Omega$ is centrally symmetric, one has
\[I_2(graph(\psi))=-graph(\psi)\, ,\]
where $\psi(t) = -L_1(t,t)$.
\item\label{punto 2} $I(\Pi) = \Pi$, where $\Pi$ is the set of $4$-periodic points whose orbit is contained in the zero section.
\item\label{punto 3} If $\Omega$ is centrally symmetric, then $I^2(\Lambda)=\Lambda$ for every $\lambda \in (0,1)$, i.e., 
$(t,s) \in \Lambda$ if and only if $(t^*,s) \in \Lambda$. 
\item\label{punto 4} If $\Omega$ is centrally symmetric, then $I(\Lambda)=\mathcal{H}_{-\frac{1}{\lambda}}(\Lambda)$ 
for every $\lambda \in (0,1)$, where $\mathcal{H}_{-\frac{1}{\lambda}}(t,s)=(t,-\frac{s}{\lambda})$.
\end{enumerate}
\end{lemma}

\begin{proof}
Point \ref{punto 1} follows from the definition of the function $I_2$, and since points of the graph of $\psi_1$ (resp. $\psi_2$) correspond to points of type $(t,t)$ (resp. $(t,t^*)$) in $\hat{\mathcal{P}}$. 
For point \ref{punto 2}, let $(t_1,0)\in\Pi$: since the origin is compatible, its orbit is given by \[\{(t_1,0),(t_2,0),(t_1^*,0),(t_2^*,0)\}\] and in particular, we have $I(t_1,0)=I(t_1,-L_1(t_1,t_2))=(t_2^*,-L_1(t_2^*,t_1))=(t_2^*,0)$.
Concerning point \ref{punto 3}, since $I$ is continuous and by the properties of the Birkhoff attractor, the set $I^2(\Lambda)$ is a compact, connected set which separates $\mathcal{P}$. Observe that, when $\Omega$ is centrally symmetric, we have that $I^2(t,s)=(t^*,s)$, since, in such a case, $L_1(t,\tilde t)=L_1(t^*,\tilde{t}^*)$. Moreover, the map $I^4$ is the identity. Therefore, one has
\[I^2\circ T_\lambda(t_1,s_1)= I^2(t_2, \lambda s_2)=(t_2^*,\lambda s_2)\]
and 
\[T_\lambda\circ I ^2(t_1,s_1)=T_\lambda(t_1^*,s_1)=(t_2^*,\lambda s_2)\, ,\]
that is, $I^2\circ T_\lambda=T_\lambda\circ I^2$. Thus, by the $T_\lambda$-invariance of $\Lambda$, we obtain $I^2(\Lambda)=I^2\circ T_\lambda(\Lambda)=T_\lambda\circ I^2(\Lambda)$, i.e., the set $I^2(\Lambda)$ is $T_\lambda$-invariant. By the minimality property of the Birkhoff attractor, we deduce that $\Lambda\subset I^2(\Lambda)$. By applying $I^2$, we conclude that $I^2(\Lambda)\subset \Lambda$, that is $\Lambda=I^2(\Lambda)$, as required.

Let us prove point \ref{punto 4}. 
We first notice that 
$\hat{I}_2=\hat{T}\circ \hat{I}$, since, by the centrally symmetric hypothesis on $\Omega$, the vector $\gamma(t_2)-\gamma(t_0)$ is parallel to the vector $\gamma(t_0^*)-\gamma(t_2^*)$, see Figure \ref{I_2}. 
\begin{figure}
    \centering
    \includegraphics[width=0.3\linewidth]{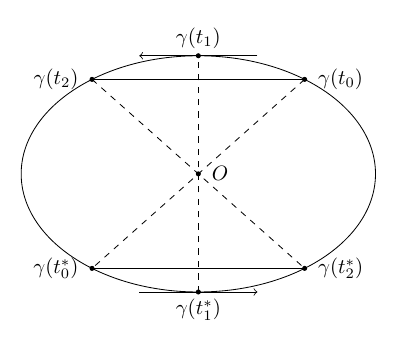}
    \caption{For the symplectic billiard map on centrally symmetric domains, $\hat{T}(t_0,t_1)=(t_1,t_2)$ if and only if $\hat{T}(t_0^*,t_1^*)=(t_1^*,t_2^*)$. In particular, $\gamma'(t_1)$ is parallel to $\gamma(t_2^*)-\gamma(t_0^*)$.}
    \label{I_2}
\end{figure}
As a consequence, 
\begin{equation} \label{prima}
I_2=T\circ I\, .
\end{equation}
In addition \begin{equation}\label{ii}
     I_2(t_1,s_1)=(t_1,-s_1)\, .
 \end{equation}
 The above equality is a consequence of the geometry and the dynamics of symplectic billiards in centrally symmetric domains, since 
 \[(t_1,s_1)=(t_1, -L_1(t_1,t_2))\stackrel{I_2}{\longmapsto}(t_1, -L_1(t_1,t_0^*))=(t_1, L_1(t_1,t_2))=(t_1,-s_1)\, ,\]
 where in the last equality, we have used the fact that
 \[L(\gamma'(t_1),\gamma(t_2)-\gamma(t_0))=L(\gamma'(t_1),\gamma(t_0^*)-\gamma(t_2^*)) = 0 \quad\implies\quad - L_1(t_1,t_0^*) = - L_1(t_1,t_2^*)= L_1(t_1,t_2)\, .\]
By equality \eqref{prima}, by point \ref{punto 3} and by the definition of $T_\lambda$, we get
\[I_2\circ I(\Lambda)=T\circ I\circ I(\Lambda)=T(\Lambda)=\mathcal{H}_{\frac{1}{\lambda}} \circ T_\lambda(\Lambda)=\mathcal{H}_{\frac{1}{\lambda}}\Lambda\, .\]
Finally, applying $I_2^{-1}$ on both sides of the previous equation and using \eqref{ii}, we conclude
\[I(\Lambda)=I_2^{-1}\biggl(\mathcal{H}_{\frac{1}{\lambda}}(\Lambda)\biggr)=\mathcal{H}_{-\frac{1}{\lambda}}(\Lambda)\, .\]
\end{proof}

\noindent We are now ready to prove the statement of point $(a)$. 
\begin{proposition} \label{punti in sezione}
Let $\Omega$ be a strongly convex, centrally symmetric $C^k$ domain, $k\geq 2$. 
Let $\lambda(\Omega)\in (0,1)$ be given by Proposition \ref{cones}. Then, for every $\lambda \in (0, \lambda(\Omega))$, one has
$\Lambda\cap(\mathbb{S}\times\{0\})= \Pi$.
\end{proposition}
\begin{proof} 
By Theorem \ref{nc}, for 
$\lambda\in(0,\lambda(\Omega))$, the Birkhoff attractor is a normally contracted graph over $\mathbb{S}$. Since it is a graph and by point \ref{punto 4} of Lemma \ref{secondo}, we deduce that $\Lambda\cap I(\Lambda)=\Lambda\cap(\mathbb{S}\times\{0\})$.
Thus, by point \ref{punto 3} of Lemma \ref{secondo}, if $(t_1,0)\in \Lambda\cap(\mathbb{S}\times\{0\})=\Lambda\cap I(\Lambda)$, then $I(t_1,0)\in I(\Lambda)\cap\Lambda=\Lambda\cap(\mathbb{S}\times\{0\})$.
By the definition of $I$, we obtain that
\[
 I(t_1,0) = I(t_1,-L_1(t_1,t_2)) = (t_2^*,-L_1(t_2^*,t_1))=(t_2^*,0)\, ;
\]
in particular, the vector $\gamma'(t_1)$ is parallel to $\gamma(t_2)-O$ and $\gamma'(t_2)$ is parallel to $\gamma(t_1)-O$. We then deduce that the point $(t_1,0)\in \Pi$, proving then $\Lambda\cap(\mathbb{S}\times\{0\})\subset \Pi$.

Moreover, again by Theorem \ref{nc}, we also know that the Birkhoff attractor coincides with the global attractor, when $\lambda\in(0,\lambda(\Omega))$. Let $(t_1,0)\in \Pi$: then clearly $(t_1,0)\in \mathbb{S}\times\{0\}$ and, since $(t_1,0)\in \Lambda_0$ and $\Lambda_0=\Lambda$, we easily conclude that $(t_1,0)\in\Lambda\cap(\mathbb{S}\times\{0\})$ as required.

\end{proof}

\subsection{An open and dense property in the centrally symmetric case}

\indent To prove point $(b)$, we first need to discuss the quantity and quality of $4$-periodic orbits for centrally symmetric domains and then to introduce the topology with respect to which we will state our result.
Let $\Omega$ be a $C^k$ centrally symmetric domain, $k\geq 2$, and let $\gamma\colon \mathbb{S}\to \mathbb{R}^2$ be a parametrization of $\partial\Omega$.
Let $(t_1,t_2)$ be a 4-periodic point. Consider then the quantity:
\begin{equation} \label{L12 theta}
k_{1,2} = k_{1,2}(t_1,t_2) := \dfrac{L_{11}(t_1,t_2)\cdot L_{22}(t_1,t_2)}{L_{12}^2(t_1,t_2)}\, .
\end{equation}
\begin{lemma}\label{bifurcation} Let $\Omega$ be a $C^k$ strictly convex centrally symmetric domain, $k\geq 2$. 
Let $(t_1,t_2)$ be a $4$-periodic point.
For $\lambda \in (0,1)$, denote by $\mu_1 = \mu_1(\lambda)$ and $\mu_2 = \mu_2(\lambda)$ the eigenvalues of $DT^4_{\lambda}(t_1,s_1(t_1,t_2))$, with $|\mu_1| \le |\mu_2|$. Then the next cases occur.  

\begin{itemize}
    \item[$(a)$] If $k_{1,2}>1$, then $0<\mu_1<\lambda^4<1<\mu_2$, and the $4$-periodic point is a \textbf{saddle}.
    \item[$(b)$] If $k_{1,2}=1$, then $\mu_1=\lambda^4$, $\mu_2=1$, and the $4$-periodic point is \textbf{parabolic}.
    \item[$(c)$] If $k_{1,2}\in (0,1)$, then the $4$-periodic point is a \textbf{sink}. In particular, let 
    $$\lambda_-=\lambda_-(t_1,t_2):=\frac{1-\sqrt{1-k_{1,2}}}{1+\sqrt{1-k_{1,2}}}\in(0,1)\, .$$
    Thus:
    \begin{itemize}
        \item[$(i)$] If $\lambda\in(0,\lambda_-)$, then $\mu_1,\mu_2\in\mathbb{R}$ and $\lambda^4<\mu_1<\mu_2<1$.
        \item[$(ii)$] If $\lambda=\lambda_-$, then $\mu_1=\mu_2=\lambda^2$.
        \item[$(iii)$] If $\lambda\in(\lambda_-,1)$ and $k_{1,2} \ne \frac{2\lambda}{(\lambda + 1 )^2}$, then $\overline{\mu_1}=\mu_2$ and $\vert \mu_1\vert=\vert \mu_2\vert=\lambda^2$. 
        \item[$(iv)$] If $\lambda\in(\lambda_-,1)$ and $k_{1,2} = \frac{2\lambda}{(\lambda + 1 )^2}$, then
        $\mu_1 = \mu_2 = - \lambda^2$. 
    \end{itemize}
\end{itemize}
\end{lemma}

\begin{proof} 
We already know by Lemma \ref{formula DT lambda} the explicit formula for the differential of the dissipative symplectic billiard map. We are studying $4$-periodic orbits, hence, up to choose a good parametrization, all of the form:
\begin{equation} \label{formula 4 per}
\{ (t_1,t_2), (t_2,t_1 + \pi), (t_1 + \pi, t_2 + \pi), (t_2 + \pi,t_1)\}\, ,
\end{equation}
see Remark \ref{tutte}. In such a case, from the central symmetric hypothesis, the next equalities hold:
$$L_{12}(t_2,t_1+\pi) = L_{12}(t_1,t_2), \quad L_{11}(t_2,t_1+\pi) = L_{22}(t_1,t_2), \quad L_{22}(t_2,t_1+\pi) = L_{11}(t_1,t_2)\,.$$
From now on, we omit the dependence of every $L_{ij}$ on $(t_1,t_2)$. By using the previous equalities, a direct computation gives:
\begin{equation*}
\begin{split}
A_{\lambda} & := DT^2_\lambda(t_1,s_1(t_1,t_2)) = DT_{\lambda} (t_2,s_1(t_2,t_1 +\pi)) \cdot DT_{\lambda} (t_1,s_1(t_1,t_2)) \\
& = \frac{1}{L_{12}^2} \begin{pmatrix}
      -\lambda L_{12}^2 + (1+\lambda) L_{11} \cdot L_{22} & (1+\lambda) L_{22} \\
(1+\lambda) (-\lambda L_{12}^2 +\lambda L_{11} \cdot L_{22}) \cdot L_{11} & -\lambda L_{12}^2 + (\lambda +\lambda^2) L_{11} \cdot L_{22}
    \end{pmatrix}\, ,
    \end{split}
\end{equation*}
and, for the $4$-periodic orbit \eqref{formula 4 per}, one has $DT^4_\lambda(t_1,s_1(t_1,t_2))=A_\lambda^2$. Consequently, to understand the nature of the $4$-periodic points of $T_\lambda$, we just need to study the eigenvalues of the matrix $A_\lambda$. The determinant and the trace of $A_\lambda$ are respectively
\begin{equation*}\label{dettr}
    \det A_\lambda=\lambda^2 \qquad \text{and} \qquad tr A_\lambda= -2 \lambda + (1 + \lambda)^2 \left( \frac{L_{11} \cdot L_{22}}{L_{12}^2} \right) \, .
\end{equation*}
The characteristic polynomial is then
 \[\chi_\lambda(x)=x^2-\left[(1 + \lambda)^2 k_{1,2} -2\lambda\right]x+\lambda^2\, ,\]
 which is exactly as in the Birkhoff dissipative case with $\lambda_1=\lambda_2=\lambda$ and $k_{1,2} := \frac{L_{11} \cdot L_{22}}{L_{12}^2}$, see \cite[Appendix A]{BFL}, where we refer to the notation of the Appendix A there. Repeating then \textit{verbatim} the proof in \cite{BFL}, we conclude.
 \end{proof}

\begin{figure}
    \centering
    \includegraphics[width=0.4\linewidth]{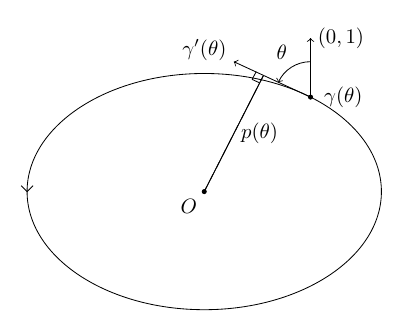}
    \caption{The support function $p(\theta)$ at the point  $\theta\in\mathbb{S}$.}
    \label{support}
\end{figure}

Let us introduce some further notations that will be largely used in the sequel. Denote by
$$e_\theta =(-\sin\theta, \cos\theta)$$
the unit vector which forms an angle $\theta\in \mathbb{S}$ with the fixed vertical direction $(0,1)$. For every $\theta$ there exists a unique point $\gamma(t_\theta):=\gamma(\theta)\in\partial\Omega$ such that $\gamma'(\theta)=\|\gamma'(\theta)\|e_\theta$. Let $O\in \mathrm{int}(\Omega)$ be the center of symmetry of the table. 
Let $p\colon \mathbb{S}\to\mathbb{R}_+$ be the support function, defined as the distance of the origin from the affine line $\gamma(\theta)+\mathbb{R}\gamma'(\theta)$. 
Let $J$ be the rotation of angle $\frac{\pi}{2}$. Then 
we have
\begin{equation}\label{diff1}
\begin{cases}
\gamma(\theta)&=p'(\theta)e_\theta-p(\theta)Je_\theta\, ;\\
    \gamma'(\theta)&=(p(\theta)+p''(\theta))e_\theta\, ;\\
    \gamma''(\theta)&=(p'(\theta)+p'''(\theta))e_\theta+(p(\theta)+p''(\theta))Je_\theta\, ;
    \end{cases}
\end{equation}
see \cite{BBN2}. 
We refer to Figure \ref{support}.
We remind that the ray of the osculating circle at $\gamma(\theta)$ is 
$$\rho(\theta) = p(\theta)+p''(\theta)\, ,$$ 
see for example \cite[Formula $2.9$]{Flan68}. We are going to use the angles $\theta$ as coordinates.
In the centrally symmetric case (see Remark \ref{tutte}) $\Pi$ is the set of all $4$-periodic points for $\{T_{\lambda}\}_{\lambda \in [0,1]}$. A point $(\theta_1,\theta_2)$ corresponding to a point in $\Pi$ has then the following 
$4$-periodic orbit:
\begin{equation} \label{Natale}
\{ (\theta_1,\theta_2), (\theta_2,\theta_1 + \pi), (\theta_1 + \pi, \theta_2 + \pi), (\theta_2 + \pi,\theta_1) \}\, .  
\end{equation}
With an abuse of notation, we will indicate the symplectic billiard map in the coordinates $(\theta_1,\theta_2)$ also by $T$.
For a 4-periodic point corresponding to the couple of angles $(\theta_1,\theta_2)$, we denote the relative quantity $k_{1,2}$, defined in \eqref{L12 theta}, also as $k_{1,2}=k_{1,2}(\theta_1,\theta_2)$. 
\begin{remark}\label{remark per formula esplicita k 1 2}
We explicit now the quantity $\frac{L_{11} L_{22}}{L_{12}^2}$, that appears in Lemma \ref{bifurcation}, in terms of $(\theta_1,\theta_2)$. We start by recalling that at the points of a $4$-periodic orbit, the vectors $\gamma'(\theta_i)$ and $\gamma(\theta_{i+1})$ are parallel, so $\mathrm{det}(e_{\theta_i},\gamma(\theta_{i+1}))=0$, that is
\begin{equation} \label{quattro periodica*}
    \begin{split}
\mathrm{det}(e_{\theta_1},\gamma(\theta_{2})) & = p'(\theta_2) \sin(\theta_2-\theta_1) - p(\theta_2) \cos(\theta_2 -\theta_1) = 0\, ,\\
\mathrm{det}(e_{\theta_2},\gamma(\theta_{1})) & = - p'(\theta_1) \sin(\theta_2-\theta_1) - p(\theta_1) \cos(\theta_2 -\theta_1) = 0\, .
    \end{split}
\end{equation}
As a consequence:
\begin{equation*} \label{quattro periodica}
    \begin{split}
L_{12}(t_{\theta_1},t_{\theta_2}) & =\rho(\theta_1)\rho(\theta_2)\sin{(\theta_2-\theta_1)}\, ,\\
L_{22}(t_{\theta_1},t_{\theta_2}) & =-p(\theta_1)\rho(\theta_2)\sin{(\theta_2-\theta_1)}+p'(\theta_1)\rho(\theta_2)\cos{(\theta_2-\theta_1)}\, ,\\
L_{11}(t_{\theta_1},t_{\theta_2}) & =-p(\theta_2)\rho(\theta_1)\sin{(\theta_2-\theta_1)}-p'(\theta_2)\rho(\theta_1)\cos{(\theta_2-\theta_1)}\, .
    \end{split}
\end{equation*}
\noindent So that
\begin{equation*} 
L_{11}(t_{\theta_1},t_{\theta_2}) \cdot L_{22}(t_{\theta_1},t_{\theta_2}) = 
\end{equation*}
\begin{equation*}
\begin{split}
\rho(\theta_1) \rho(\theta_2) \left[ p(\theta_1) p(\theta_2) \sin^2(\theta_2 -\theta_1) - p'(\theta_1) p'(\theta_2) \cos^2 (\theta_2 - \theta_1)  +\right.\\
\left.\left(p(\theta_1)p'(\theta_2)-  p'(\theta_1) p(\theta_2)\right) \sin(\theta_2 -\theta_1) \cos(\theta_2 -\theta_1) \right] \, .
\end{split}
\end{equation*}
Using now the equalities \eqref{quattro periodica*}, previous formula simplifies to 
\begin{equation} \label{L11 L22}
 L_{11}(t_{\theta_1},t_{\theta_2}) \cdot L_{22}(t_{\theta_1},t_{\theta_2}) = \rho(\theta_1) \rho(\theta_2) \left[ p(\theta_1) p(\theta_2) -p'(\theta_1) p'(\theta_2) \right] = \frac{\rho(\theta_1) \rho(\theta_2) p(\theta_1) p(\theta_2)}{\sin^2(\theta_2 -\theta_1)}\, .
\end{equation}
Moreover:
\begin{equation} \label{L12}
L_{12}^2(t_{\theta_1},t_{\theta_2}) = \rho^2(\theta_1) \rho^2(\theta_2) \sin^2(\theta_2 - \theta_1)\, .
\end{equation}
From \eqref{L11 L22} and \eqref{L12}, we have that
$$\frac{L_{11} L_{22}}{L_{12}^2} (\theta_1,\theta_2)= \frac{p(\theta_1) p(\theta_2)}{\rho(\theta_1) \rho(\theta_2) \sin^4 (\theta_2 -\theta_1)}\, .$$
By using again equalities \eqref{quattro periodica*}, we obtain:
$$(p')^2(\theta_i) \sin^2(\theta_2-\theta_1) = p^2(\theta_i) \cos^2(\theta_2 -\theta_1) = p^2(\theta_i) (1 - \sin^2(\theta_2 -\theta_1)) \Rightarrow \sin^2(\theta_2-\theta_1) = \frac{p^2(\theta_i)}{((p')^2 + p^2)(\theta_i)}\, ,$$
for $i = 1,2$. We then conclude that 
\begin{equation}\label{k 12 con angoli} 
 k_{1,2} := \frac{L_{11} L_{22}}{L_{12}^2} = \frac{(p'^2 + p^2)(\theta_1) (p'^2 + p^2)(\theta_2)}{{\rho(\theta_1)\rho(\theta_2)} p(\theta_1) p(\theta_2)}\, .
\end{equation}
\end{remark} 



In the following lemma, we see how it is possible to perturb any strongly convex, centrally symmetric domain in order to assure that all the 4-periodic points of the dynamics of the perturbed table are non-degenerate (in particular, they are only a finite number). Recall that a 4-periodic point is non-degenerate for the map $T$, if the differential map $DT^4$ at the point does not admit the values $\pm 1$ as eigenvalues.

\begin{lemma} \label{lemma1}
Let $\Omega$ 
be a $C^k$ strongly convex, centrally symmetric domain, $k\geq 2$, and let \( p\colon \mathbb{S}\rightarrow \mathbb{R} \) be its support function of class \( C^k \). Then, for every \( \varepsilon > 0 \), there exists a strongly convex, centrally symmetric smooth domain $\Omega_\varepsilon$ with support function \( p_\varepsilon\colon \mathbb{S}\to \mathbb{R} \) such that
\[
\|p - p_\varepsilon\|_{k} < \varepsilon\, ,
\]
where $\|\cdot \|_{k}$ denotes the $C^k$-norm, and such that the symplectic billiard associated with the domain \( \Omega_\varepsilon \)  has a finite number of 4-periodic orbits, all of which are non-degenerate.
\end{lemma}

\begin{proof}
Let $(\theta_1,\theta_2)$ correspond to a $4$-periodic point in $\Pi$; then its orbit is of the
form \eqref{Natale}. In particular, the vector $\gamma'(\theta_{2})$ is parallel to the vector $\gamma(\theta_1) - O$, and so, by \eqref{diff1}, one has 
$\det(e_{\theta_2}, \gamma(\theta_1)) = 0$, that is 
\[p(\theta_1) \cos(\theta_2 - \theta_1) + p'(\theta_1) \sin(\theta_2 - \theta_1) = 0 \quad \Leftrightarrow \quad \tan\left((\theta_2-\theta_1)-\frac{\pi}{2}\right)=\frac{p'}{p}(\theta_1)\, ,\]
where we use the fact that the support function $p$ is positive and that $\theta_2\neq \theta_1$. There exist lifts 
of the angles $\theta_1,\theta_2$ --which, for the sake of simplicity, we continue to indicate by $\theta_1,\theta_2$-- such that $\theta_2 - \theta_1\in(0,\pi)$. Therefore, we have
\begin{equation} \label{secondo angolo}
\theta_2=\arctan\left(\frac{p'}{p}(\theta_1)\right)+\frac{\pi}{2}+\theta_1\, .
\end{equation}
Apply then the same argument to the subsequent points of the orbit: since $\theta_3=\theta_1+\pi$, we obtain that the angle $\theta_1$ has to solve
    \begin{equation} \label{4per}
        G(\theta) := \dfrac{p'}{p} \left( \arctan \left(\frac{p'}{p} (\theta) \right) + \theta +\frac{\pi}{2} \right) + \dfrac{p'}{p}(\theta)=0\, .
    \end{equation}
\indent We proceed by observing that, up to perturb the support function (and so the billiard table), equation \eqref{4per} admits a finite number of solutions. Let $(p_n)_{n\in\mathbb{N}}$ be a sequence of $\pi$-periodic trigonometric polynomials that approximates the support function $p$ in the $C^k$-norm. Substitute then $p_n$ in \eqref{4per}: one obtains an equation whose left-hand side is a real analytic function, which therefore has a finite number of solutions. This fact immediately gives that every strictly convex domain can be perturbed into a strictly convex domain with a finite number of 4-periodic orbits.  \\
\indent The last part of the proof is devoted to prove that it is possible to further perturb the domain in such a way that all $4$-periodic orbits become non-degenerate. Let $n\in\mathbb{N}$ be large enough and fix the support function $p_n$. With an abuse of notation, we still denote by $\Omega$ the billiard table corresponding to the support function $p_n$ and we designate $p_n$ by $p$.
Let \( g(\theta) \) be a $\pi$-periodic function to be chosen later. We consider the perturbed domain whose support function is given by \( p_\varepsilon(\theta):=e^{\varepsilon g(\theta)} p(\theta) \), and define the function
\[
f(\theta,\varepsilon):= \frac{p'}{p}(\theta)+\varepsilon g'(\theta)\, .
\]
Observe that $\frac{p'_\varepsilon}{p_\varepsilon} (\theta) =f(\theta,\varepsilon)$.
The 4-periodic orbits of the billiard table associated to the support function $p_\varepsilon$ are then determined by solving the following equation, analogous to\eqref{4per}:
\begin{equation} \label{pert1}
    G(\theta,\varepsilon):=f\left(\arctan(f(\theta,\varepsilon))+\theta+\frac{\pi}{2},\varepsilon\right)+f(\theta,\varepsilon)=0\, ,
\end{equation}
that we can write also as
\begin{equation}\label{pert1 v2}
G(\theta,\varepsilon)=\dfrac{p'_\varepsilon}{p_\varepsilon}\left(\arctan\left(\dfrac{p'_\varepsilon}{p_\varepsilon}(\theta)\right)+\theta+\frac{\pi}{2}\right)+\dfrac{p'_\varepsilon}{p_\varepsilon}(\theta)=0\, .
\end{equation}
Let us show that a 4-periodic orbit for the table with support function $p_\varepsilon$ is degenerate if and only if $\partial_\theta G(\theta,\varepsilon)=0$. Let $\theta_1$ correspond to a 4-periodic point for the table associated to the support function $p_\varepsilon$, i.e., $G(\theta_1,\varepsilon)=0$. The point is degenerate if and only if, according to formula \eqref{k 12 con angoli} and point $(b)$ of Lemma \ref{bifurcation}, we have
\begin{equation}\label{k 1 2 avec p epsilon}
k_{1,2}=\frac{(p'^2_\varepsilon + p^2_\varepsilon)(\theta_1) (p'^2_\varepsilon + p^2_\varepsilon)(\theta_2)}{{\rho_\varepsilon(\theta_1)\rho_\varepsilon(\theta_2)} p_\varepsilon(\theta_1) p_\varepsilon(\theta_2)}=1\, ,
\end{equation}
where $\rho_\varepsilon=p_\varepsilon+p''_\varepsilon$ and, with an abuse of notation, we still denote by $\theta_2$ the quantity $\arctan\left(\frac{p'_\varepsilon}{p_\varepsilon}(\theta_1)\right)+\theta_1+\frac{\pi}{2}$. Consider then
\begin{equation*}
\begin{split}
\partial_\theta G(\theta_1,\varepsilon)&= \frac{p''_\varepsilon p_\varepsilon - p'^2_\varepsilon}{p^2_\varepsilon}(\theta_2) \left( \frac{1}{1 + \frac{p'^2_\varepsilon}{p^2_\varepsilon} (\theta_1)} \left(\frac{p''_\varepsilon p_\varepsilon-p'^2_\varepsilon}{p^2_\varepsilon}\right)(\theta_1) + 1 \right) + \frac{p''_\varepsilon p_\varepsilon-p'^2_\varepsilon}{p^2_\varepsilon}(\theta_1)\\
&= \frac{p''_\varepsilon p_\varepsilon - p'^2_\varepsilon}{p^2_\varepsilon}(\theta_2) \left( \frac{(p''_\varepsilon + p_\varepsilon) p_\varepsilon}{p'^2_\varepsilon + p^2_\varepsilon}(\theta_1) \right) + \frac{p''_\varepsilon p_\varepsilon-p'^2_\varepsilon}{p^2_\varepsilon}(\theta_1) \\
& = \frac{(p''_\varepsilon + p_\varepsilon) p_\varepsilon - (p'^2_\varepsilon + p^2_\varepsilon)}{p^2_\varepsilon}(\theta_2) \left( \frac{(p''_\varepsilon + p_\varepsilon) p_\varepsilon}{p'^2_\varepsilon + p^2_\varepsilon}(\theta_1) \right) + \frac{(p''_\varepsilon + p_\varepsilon) p_\varepsilon- (p'^2_\varepsilon + p^2_\varepsilon)}{p^2_\varepsilon}(\theta_1) \\
& = \left( \frac{\rho_\varepsilon}{p_\varepsilon}(\theta_2) - \frac{p'^2_\varepsilon + p^2_\varepsilon}{p^2_\varepsilon}(\theta_2) \right) \left( \frac{\rho_\varepsilon p_\varepsilon}{p'^2_\varepsilon + p^2_\varepsilon}(\theta_1) \right) + \frac{\rho_\varepsilon}{p_\varepsilon}(\theta_1) - \frac{p'^2_\varepsilon + p^2_\varepsilon}{p^2_\varepsilon}(\theta_1)\, .
\end{split}
\end{equation*}
\noindent From \eqref{k 1 2 avec p epsilon}, we have that $\frac{p'^2_\varepsilon+p^2_\varepsilon}{\rho_\varepsilon p_\varepsilon}(\theta_2)=\frac{\rho_\varepsilon p_\varepsilon}{p'^2_\varepsilon+p^2_\varepsilon}(\theta_1)$, and we obtain from the previous equalities:{\color{red}ANNA QUI}
\begin{equation*}
\begin{split}
\partial_\theta G(\theta_1,\varepsilon)&= \dfrac{p'^2_\varepsilon+p^2_\varepsilon}{p^2_\varepsilon}(\theta_2)\left[ 1-\dfrac{\rho_\varepsilon p_\varepsilon}{p'^2_\varepsilon+p^2_\varepsilon}(\theta_1) \right]-\dfrac{p'^2_\varepsilon+p^2_\varepsilon}{p^2_\varepsilon}(\theta_1)\left[1-\dfrac{\rho_\varepsilon p_\varepsilon}{p'^2_\varepsilon+p^2_\varepsilon}(\theta_1)\right]\\
&=\left[1-\dfrac{\rho_\varepsilon p_\varepsilon}{p'^2_\varepsilon+p^2_\varepsilon}(\theta_1)\right]\, \left[\dfrac{p'^2_\varepsilon+p^2_\varepsilon}{p^2_\varepsilon}(\theta_2)-\dfrac{p'^2_\varepsilon+p^2_\varepsilon}{p^2_\varepsilon}(\theta_1)\right]\\
&= \left[1-\dfrac{\rho_\varepsilon p_\varepsilon}{p'^2_\varepsilon+p^2_\varepsilon}(\theta_1)\right]\,\left[ \dfrac{p'^2_\varepsilon}{p^2_\varepsilon}(\theta_2)-\dfrac{p'^2_\varepsilon}{p^2_\varepsilon}(\theta_1)\right]\, .
\end{split}
\end{equation*}
Since we are looking at a point satisfying \eqref{pert1 v2}, it holds $\frac{p'_\varepsilon}{p_\varepsilon}(\theta_2)=-\frac{p'_\varepsilon}{p_\varepsilon}(\theta_1)$, and we conclude that
\begin{equation}\label{pert2}
\partial_{\theta}G(\theta_1,\varepsilon)=0\, ,
\end{equation}
for every $\theta_1$ that corresponds to a degenerate 4-periodic point.
Note that both conditions, \eqref{pert1 v2} and \eqref{pert2}, are expressed by real-analytic equations. \\
If no solution of \eqref{pert1} also satisfies \eqref{pert2} for sufficiently small \( \varepsilon \), we are done. Suppose instead, by contradiction, that this is not the case, i.e., there exists a sequence \( \varepsilon_n \to 0 \) and corresponding values \( \theta_n \to \theta_0 \) such that $G(\theta_n,\varepsilon_n)=\partial_\theta G(\theta_n,\varepsilon_n)=0$ for every $n$.
Clearly, \( \theta_0 \) corresponds to a degenerate 4-periodic orbit of the unperturbed domain \( \Omega \). Note that
\[
\partial_\varepsilon G(\theta,0)
= g'\left(\arctan(f(\theta,0))+\theta+\frac{\pi}{2}\right)
+ g'(\theta) \left( \frac{\partial_\theta f\left(\arctan(f(\theta,0))+\theta+\frac{\pi}{2},0\right)}{1+f^2(\theta,0)} + 1 \right)\, .
\]
Choose \( g \) so that \( \partial_\varepsilon G(\theta_0,0)\neq 0 \), for every $\theta_0$ corresponding to a 4-periodic orbit of the table $\Omega$.
This is possible because we only need to check a finite number of points.

Thanks to the assumptions on \( g \), by the implicit function theorem we have that the solutions of \eqref{pert1} in a neighborhood of \( \theta_0 \) are given by the graph of a real-analytic function \( \varepsilon(\theta) \). Since $\Omega$ has only a finite number of 4-periodic orbits, the function \( \varepsilon(\theta) \) is not identically zero. There exists then \( N \in \mathbb{N} \) and a function $c$, which does not vanish in a neighborhood of any $\theta_0$, corresponding to a 4-periodic orbit for $\Omega$, such that
\begin{equation} \label{pert3}
    \varepsilon(\theta) = (\theta - \theta_0)^N c(\theta)\, .
\end{equation}
We then have
\[
\varepsilon_n = c(\theta_n)(\theta_n - \theta_0)^N\, ,
\]
and, by assumption,
\[
\partial_\theta G(\theta_n,\varepsilon(\theta_n)) = 0\, .
\]
It follows, by real-analyticity, that
\[
\partial_\theta G(\theta, \varepsilon(\theta)) = 0 \quad \text{for all } \theta\, .
\]
Then, again by the implicit function theorem, in a neighborhood of the graph of \( \varepsilon(\cdot) \) we can write
\[
\partial_\theta G(\theta, \varepsilon) = G(\theta, \varepsilon)\, h(\theta, \varepsilon)\, ,
\]
for some suitable function \( h \). Differentiating with respect to \( \theta \), we get
\[
\partial^2_\theta G(\theta, \varepsilon) = G(\theta, \varepsilon) \left(h^2(\theta, \varepsilon) + \partial_\theta h(\theta, \varepsilon)\right)\, ,
\]
and proceeding inductively, we find
\[
\partial^n_\theta G(\theta, \varepsilon(\theta)) = 0 \quad \text{for all } n\, .
\]
As a consequence, by real-analyticity, we deduce that \( G \) is independent of \( \theta \) in an open set by \eqref{pert3}, and hence everywhere. This is the required contradiction and concludes the proof.
\end{proof}
\noindent We want now to make precise the topology used in order to be able to talk about open and dense sets in point $(b)$. Denote by $C^\infty(\mathbb{S},\mathbb{R}^2)$ the set of smooth functions $\gamma\colon\mathbb{S}\to \mathbb{R}^2$ and endow it with the $C^2$-norm $\Vert\cdot\Vert_2$. Consider the open set 
\[
\mathcal{B}^0_2:=\{\gamma\in C^\infty(\mathbb{S},\mathbb{R}^2) :\ \gamma\colon\mathbb{S}\to \mathbb{R}^2\text{ is an embedding}\}\, ;
\]
we will be interested in the subset 
\[
\mathcal{B}_2:=\{\gamma\in\mathcal{B}^0_2 :\ \gamma(\mathbb{S})\text{ is strongly convex and centrally symmetric}\}
\]
endowed with the restricted topology induced by the norm $\Vert\cdot\Vert_2$. 
The following result is then immediate from Lemma \ref{lemma1}. 
\begin{corollary}\label{corollary finite non deg 4 per}
    There exists an open and dense set $\mathcal{U}\subset \mathcal{B}_2$ such that, for every $\gamma\in\mathcal{U}$, the symplectic billiard map associated to $\gamma$ has a finite number of 4-periodic orbits, all non-degenerate.
\end{corollary}

We are now ready to prove the following statement, which corresponds to point $(b)$.
\begin{proposition} \label{unione varieta instabili}
There exists an open and dense set $\mathcal{U}\subset \mathcal{B}_2$ such that for any $\gamma\in\mathcal{U}$ the following property holds. Let $\lambda'(\Omega) \in (0, 1)$ as given in Theorem \ref{nc}; there exists $\lambda''(\Omega)\in(0,\lambda'(\Omega))$ such that for every $\lambda\in(0,\lambda''(\Omega))$, the Birkhoff attractor $\Lambda$ has rotation number $1/4$ and 
\[\Lambda=\bigcup_{i=1}^l \bigcup_{j=0}^3\overline{\mathcal{W}^u(T_\lambda^j(H_i);T_\lambda^4)}\]
for some finite collection $\lbrace H_i\rbrace_{i=1}^l$ of $4$-periodic points of saddle type, where $\mathcal{W}^u(H_i;T_\lambda^4)$ denotes the unstable manifold of a hyperbolic point $H_i$ with respect to the dynamics $T^4_\lambda$.
\end{proposition}
\begin{proof} The proof is an adaptation of the proof of Theorem $5.14$ in \cite{BFL}. Let $\lambda'(\Omega)\in(0,1)$ be given by Theorem \ref{nc}. Then, for every $\lambda\in(0,\lambda'(\Omega))$, the (Birkhoff) attractor is a $C^1$ normally contracted graph of a function $\eta_\lambda: \mathbb{S} \to \mathbb{R}$. Let us define 
       $$ g_\lambda:\mathbb{S}\to\mathbb{S}, \qquad
            t_1\mapsto p_1\circ T_\lambda(t_1,\eta_\lambda(t_1))$$ 
     the circle map induced by $T_\lambda$ on the attractor, where $p_1$ is the projection on the first coordinate. 
     As proved in Claim $5.15$ in \cite{BFL}, 
$g_\lambda\to g_0$ in the $C^1$ topology when $\lambda\to0$. We proceed to prove that $g_0$ is a circle diffeomorphism. By Lemma \ref{formula DT lambda}, we get
     \begin{equation*}
         g'_\lambda(t_1)=-\frac{1}{L_{12}(t_1,t_2)}\left[L_{11}(t_1,t_2)+\eta'_\lambda(t_1)\right]\, .
     \end{equation*}
     In particular $g_0'(t_1)=-\frac{L_{11}(t_1,t_2)}{L_{12}(t_1,t_2)}$. 
     
We want to show that $g_0$ is a circle diffeomorphism: in particular, it is sufficient to assure that $g_0'\neq 0$ at every $t_1$. 
Its expression is, more explicitely,
\begin{equation}\label{g' no equiaff}
g'_0(t_1)=\dfrac{p}{\rho}(\theta_2)+\dfrac{p'}{\rho}(\theta_2)\cot{(\theta_2-\theta_1)}\, .
\end{equation}
 We recall that $\Pi \subset \mathbb{S} \times \{0\}$ does not depend on $\lambda \in [0,1]$. Moreover, by Lemma \ref{lemma1}, for a $C^k$-open and dense subset of centrally symmetric billiard tables, $\Pi$ consists of non-degenerate $4$-periodic points, i.e., saddles or sinks, hence persisting under perturbations. In the sequel, in order to conclude, we prove that the circle dynamics of $g_{\lambda}$ (which is a small perturbation of $g_0$ in the $C^1$ topology) essentially recovers the one of $T_{\lambda}$ on $\Lambda$. Let $(t_1,0) \in \Pi$. In particular, we have $\cot(\theta_2-\theta_1)=\frac{p'}{p}(\theta_2)$ and $\cot(\theta_1-\theta_2)=\frac{p'}{p}(\theta_1)$, thanks to the centrally symmetric hypothesis of the table. Then, according to \eqref{g' no equiaff}:
 $$(g_0^4)'(t_1) = (g'_0(t_1) g'_0(t_2))^ 2 = \left(\dfrac{p'^2+p^2}{\rho p}(\theta_2)\, \dfrac{p'^2+p^2}{\rho p}(\theta_1)\right)^2=k_{1,2}^2 \ne 1$$ 
 so that, for the circle diffeomorphism $g_0$, the $4$-periodic point $t_1$ is repelling or attracting (respectively when $k_{1,2}>1$ or $k_{1,2}<1$). Since for $\lambda>0$ small enough $g_\lambda$ is $C^1$-close to $g_0$ and the points of $\Pi$ are generically isolated, there exists $\lambda''(\Omega)\in(0,\lambda'(\Omega))$ such that for every $\lambda\in(0,\lambda''(\Omega))$ the $4$-periodic point $t_1$ for $g_0$ admits a continuation for $g_\lambda$. Thus, the restriction of $T_\lambda$ on $\Lambda$ still has rotation number $1/4$. As in the end of the proof of Theorem $5.14$ in \cite{BFL}, by considering the $\alpha$-limit and the $\omega$-limit sets of the points in $\Lambda\setminus \Pi$, we conclude that $$\Lambda_\lambda=\bigcup_{i=1}^l \bigcup_{j=0}^3\overline{\mathcal{W}^u(T_\lambda^j(H_i);T_\lambda^4)}$$
   for some finite collection $\lbrace H_i\rbrace_{i=1}^l$ of $4$-periodic points of saddle type.
\end{proof}

\subsection{Fragility of invariant curves of rotation number \texorpdfstring{$1/4$}{TEXT}}
In this subsection, we will see that, among centrally symmetric tables, invariant curves of rotation number $\frac 1 4$, for (conservative) symplectic billiard maps, are very easy to destroy. This result is used in Section \ref{birkhoff-section-complicated} to obtain topologically complicated Birkhoff attractors in the dissipative framework.

\begin{proposition}\label{open and dense no 1 4 inv curve}
There exists an open and dense set of strongly convex, centrally symmetric billiard tables $\mathcal{U}\subset \mathcal{B}_2$ such that, for every $\gamma\in\mathcal{U}$, the associated symplectic billiard map does not have an invariant curve of rotation number $\frac 1 4$.
\end{proposition}

We split the proof of Proposition \ref{open and dense no 1 4 inv curve} in two parts, in order to show first the openness and then the density property.

\begin{lemma}\label{open no inv curve}
The set of strongly convex, centrally symmetric billiard tables in $\mathcal{B}_2$ whose associated symplectic billiard map has an invariant curve of rotation number $\frac 1 4$ is closed.
\end{lemma}
\begin{proof}
Let $(\gamma_n)_{n\in\mathbb{N}}$ be a sequence of centrally symmetric tables and denote by $T_n$ the symplectic billiard map associated to $\gamma_n$. Assume that every $T_n$ exhibits an invariant curve $\Gamma_n$ of rotation number $\frac 1 4$. Moreover, assume that the sequence of $\gamma_n$ is converging to $\gamma$ in the $C^2$ topology. The sequence of corresponding billiard maps $(T_n)_{n\in\mathbb{N}}$ is going to the billiard map $T$, associated to $\gamma$, in the $C^1$ topology.
By Birkhoff's theorem \cite{B22}, since each $T_n$ is a twist map, each invariant curve $\Gamma_n$ is a graph of a Lipschitz function over $\mathbb{S}$. Moreover, up to consider $n$ large enough, since the Lipschitz constant depends on the twist condition and since the tables are converging in the $C^2$ topology, there exists a constant $C>0$ such that all invariant curves $\Gamma_n$ are graphs of $C$-Lipschitz functions. Therefore, up to subsequences, $\Gamma_n$ converges to a curve $\Gamma$, which is the graph of a $C$-Lipschitz function and is $T$-invariant. Moreover, the curve $\Gamma$ has still rotation number $\frac 1 4$. We deduce that the set of tables with invariant curves of rotation number $\frac 1 4$ is closed.
\end{proof}

\begin{proposition}\label{dense no inv curve}
The set of strongly convex, centrally symmetric billiard tables whose associated symplectic billiard map does not have an invariant curve of rotation number $\frac 1 4$ is dense among $\mathcal{B}_2$.
\end{proposition}
\begin{proof}
Let $\Omega$ be a centrally symmetric, strictly convex smooth table with an invariant curve of rotation number $1/4$. By Corollary \ref{corollary finite non deg 4 per}, we can perform a first $C^2$ perturbation to obtain a table $\Omega'$ which has a finite number of 4-periodic orbits, all non-degenerate. If the billiard map associated to $\Omega'$ has no invariant curve of rotation number $1/4$, then we are done. Assume that this is not the case, i.e., the billiard map associated to $\Omega'$ has an invariant curve $\Gamma$ of rotation number $1/4$: in particular, there are a finite number of 4-periodic points on $\Gamma$. Since $\Gamma$ is invariant, we deduce that all the 4-periodic points on the curve are hyperbolic ones. Indeed, since the points are non-degenerate, they are either hyperbolic (i.e., the corresponding differential map $DT^4$ has eigenvalues of modulus different from 1) or elliptic ones with irrational angle (i.e., the eigenvalues of $DT^4$ are $e^{2 \pi i\rho}$ with $\rho \in \mathbb{R} \setminus \mathbb{Q}$). However, if by contradiction a point is elliptic, then, after some iterations, the vertical vector $(0,1)$ would be sent by the differential dynamics to a vector pointing downward: this would contradict the existence of the invariant curve and the orientation preserving hypothesis of the map. By Poincaré classification theorem for circle homeomorphisms, the curve $\Gamma$ is composed by heteroclinic or homoclinic connections between 4-periodic hyperbolic points, see Figure \ref{InvCurNotAllPP}.
\begin{figure}[ht]
    \centering
    \includegraphics[scale=0.7]{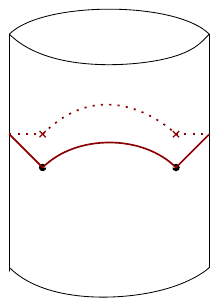}
    \caption{Invariant curve of rotation number $\frac 1 4$ with homoclinic (heteroclinic) connexions.}
    \label{InvCurNotAllPP}
\end{figure}

Denote by $T$ the symplectic billiard map associated to the table. With an abuse of notation, denote by $T$ also a lift of it. According to the notation used in \cite{Bang}, we can consider bi-infinite configurations $(t_i)_{i\in\mathbb{Z}}\in \mathbb{R}^\mathbb{Z}$. Orbits for the map $T$ are then in correspondence, thanks to the twist condition, to configurations that are critical for the formal action functional
$\mathcal{L}\colon \mathbb{R}^\mathbb{Z} \to \mathbb{R}, \ (t_i)_{i\in\mathbb{Z}} \mapsto \sum_{i\in\mathbb{Z}}L(\gamma(t_i),\gamma(t_{i+1}))$.
That is, a bi-infinite sequence corresponds to an orbit if and only if
\begin{equation}\label{caratterizzazione orbite}
L_1(t_i,t_{i+1})+L_2(t_{i-1},t_i)=0\qquad\forall i\in\mathbb{Z}\, .
\end{equation}
Let $\gamma\colon \mathbb{S}\to \mathbb{R}^2$ be the smooth arc-length parametrization of $\Omega'$, $|\partial \Omega'| = 2\pi$. Fix $\varepsilon>0$ and consider the perturbed domain corresponding to
\[
\bar\gamma\colon t\in\mathbb{S}\to\bar\gamma(t):=\gamma(t)+\delta(t)\, n(t)\in\mathbb{R}^2\, ,
\]
where $n(t)$ is the unit normal vector at $\gamma(t)$, pointing inside the domain $\Omega'$. Up to choose the smooth function $\delta$ small enough, we obtain a perturbation of the initial table which is $\varepsilon$-$C^2$ close to $\Omega'$. Our aim is then proving that we can choose the function $\delta$ so that the map of the perturbed domain does not have an invariant curve or rotation number $1/4$.

Consider $(t^-_i)_{i\in\mathbb{Z}},(t^+_i)_{i\in\mathbb{Z}}$ two sequences corresponding to nearby 4-periodic hyperbolic points contained in the invariant curve $\Gamma$. The curve $\Gamma$ is the graph of a Lipschitz function $\eta\colon\mathbb{S}\to\mathbb{R}$. Let $(t_i)_{i\in\mathbb{Z}}$ be a sequence corresponding to a point $P_0 = (t_0,\eta(t_0))\in \Gamma$, which belongs to the unstable manifold of $(t^-_0,\eta(t^-_0)) = (t^-_0,0)$ and to the stable manifold of $(t^+_0,\eta(t^+_0)) = (t^+_0,0)$. In particular it holds $t^-_i<t_i<t^+_i<t^-_i+\pi$ for every $i\in\mathbb{Z}$ and
\[
\lim_{i\to +\infty}t^+_i-t_i=\lim_{i\to -\infty}t_i-t^-_i=0\, .
\]
Let $U\subset \mathbb{R}$ be a neighborhood of $t_0$ such that $t_i\notin U$ for every $i\neq 0$. This is possible because, at $\pm\infty$, the $t_i$ accumulate on 4-periodic points, which are finite. Assume that the function $\delta$ is a positive $\pi$-periodic function whose support is contained in $U$ (i.e., its projection on $\mathbb{S}$) and such that $\delta(t_0)=\delta'(t_0)=0$ and $\delta''(t_0)\neq 0$.

Observe in particular that all the 4-periodic points, which were maximizing the area of the inscribed quadrilateral, are still 4-periodic orbits of the perturbed symplectic billiard map and moreover they are still maximizing the area of the inscribed quadrilateral for the perturbed domain, i.e., they still belong to the Aubry-Mather set of rotation number $1/4$.

Since the characterization of orbits \eqref{caratterizzazione orbite} only involves $\gamma$ and $\gamma'$ along the bi-infinite sequence corresponding to the orbit, and since both $\bar\gamma$ and $\bar\gamma'$ do not change along the sequence $(t_i)_{i\in\mathbb{Z}}$ by hypothesis, we deduce that the sequence $(t_i)_{i\in\mathbb{Z}}$ corresponds to an orbit also for the perturbed symplectic billiard map. Moreover, it is still a heteroclinc or homoclinic orbit, intersection of the unstable manifold of $(t^-_0,0)$ and the stable manifold of $(t^+_0,0)$ (for the perturbed dynamics). 

Therefore, we only have to show that this intersection is transverse. Observe that, for $\Omega'$, the tangent spaces to the stable and unstable manifolds at $(t_0,\eta(t_0))$ clearly coincide and are generated by the vector $(1,\eta'(t_0))$.

From now on, let $P_{i} := (t_i,\eta(t_i))$. Consider the perturbed dynamics and the point $P_1$: the tangent space to the stable manifold of $(t^+_1,\eta(t^+_1))$ at the point $P_1$ is still generated by the vector $(1,\eta'(t_1))$, since the dynamics remains unchanged on the future of the point.

Denote then by $T_\gamma$ and $T_{\bar\gamma}$ the symplectic billiard maps associated to $\gamma$ and $\bar\gamma$ respectively. It is then sufficient to show that
$$T_{P_1} W^s_{\bar{\gamma}}(t_1^+,0) \pitchfork T_{P_1} W^u_{\bar{\gamma}}(t_1^-,0) \,.$$
As explained right above:
$$\langle w \rangle := \langle (1,\eta'(t_1))\rangle = T_{P_1} W^s_{\bar{\gamma}}(t_1^+,0) = T_{P_1} W^s_{\gamma}(t_1^+,0)\, .$$
Moreover, let denote 
$$\langle u \rangle := T_{P_1} W^u_{\bar{\gamma}}(t_1^-,0)\, .$$ 
Then we have:
\[ w=
DT_\gamma(t_{0},s_0)\, DT_\gamma(t_{-1},s_{-1})\, \begin{pmatrix}
    1 \\ \eta'(t_{-1})
\end{pmatrix}\in \langle w\rangle \] 
and
\[u=
DT_{\bar\gamma}(t_{0},s_0)\, DT_{\bar\gamma}(t_{-1},s_{-1})\, \begin{pmatrix}
    1 \\ \eta'(t_{-1})
\end{pmatrix}\in \langle u\rangle\, .
\]
We need to prove that $w$ and $u$ are linearly independent. We start by recalling that for $i \in \mathbb{Z}$ (see formula (\ref{differential}) with $\lambda = 1$):
$$DT_\gamma (t_i,s_i) = -\frac{1}{L_{12}(t_i,t_{i+1})}\begin{pmatrix}
      L_{11}(t_i,t_{i+1}) & 1\\
- L_{12}^2(t_i,t_{i+1})+ L_{22}(t_i,t_{i+1}) L_{11}(t_i,t_{i+1})& L_{22}(t_i,t_{i+1})
    \end{pmatrix}\, .$$
Consequently, by denoting 
$$A_\gamma := L_{12}(t_{-1},t_{0}) L_{12}(t_0,t_{1}) DT_\gamma(t_{0},s_0)\, DT_\gamma(t_{-1},s_{-1})\, ,$$ 
we have that:
$$A_\gamma = \begin{pmatrix}
      L_{11}(t_{-1},t_0) L_{11}(t_0,t_{1}) + F(t_{-1},t_0) & L_{22}(t_{-1},t_0) + L_{11}(t_0,t_1)  \\
L_{11}(t_{-1},t_0)F(t_0,t_1) + L_{22}(t_0,t_1) F(t_{-1},t_0) & L_{22}(t_{-1},t_0) L_{22}(t_0,t_{1}) + F(t_{0},t_1)
    \end{pmatrix}\, ,$$
where, for $i \in \mathbb{Z}$:
$$F(t_i,t_{i+1}) := - L_{12}^2(t_i,t_{i+1})+ L_{22}(t_i,t_{i+1}) L_{11}(t_i,t_{i+1})\, .$$
Similarly, let denote 
$$B_{\bar\gamma} := L_{12}(t_{-1},t_{0}) L_{12}(t_0,t_{1}) DT_{\bar \gamma}(t_{0},s_0)\, DT_{\bar \gamma}(t_{-1},s_{-1})\, .$$ 
Observe that, since we have choosen the arc-length parametrization and $\delta(t_0)=\delta'(t_0)=0$, we get
$$
\bar{\gamma}''(t_0)=\gamma''(t_0)+\delta''(t_0)n(t_0)\, .
$$
By straightforward computations, we obtain that
\begin{equation} \label{B gamma}
B_{\bar \gamma} = A_\gamma + \det (\delta''(t_0) n(t_0), \gamma(t_1) - \gamma(t_{-1}))  C_{\gamma}  
\end{equation}
where 
$$C_{\gamma} = \begin{pmatrix}
      L_{11}(t_{-1},t_0)  & 1  \\
L_{11}(t_{-1},t_0) L_{22}(t_0,t_1) & L_{22}(t_0,t_{1})
    \end{pmatrix}\, . $$
Moreover, we notice that also $A_{\gamma}$ can be written in terms of $C_{\gamma}$:
\begin{equation} \label{A gamma}
A_{\gamma} = (L_{22}(t_{-1},t_0) + L_{11}(t_0,t_1)) C_{\gamma} + D_{\gamma}\, ,   
\end{equation}
where
$$D_{\gamma} = \begin{pmatrix}
      -L_{12}^2(t_{-1},t_{0}) & 0\\
- L_{12}^2(t_{0},t_{1}) L_{11}(t_{-1},t_0) - L_{12}^2(t_{-1},t_{0}) L_{22}(t_0,t_{1}) & -L_{12}^2(t_0,t_{1})
    \end{pmatrix}\, .$$
As a consequence of (\ref{A gamma}) and (\ref{B gamma}), $w \in \langle A_\gamma (1 , \eta'(t_{-1})) \rangle $ and $u \in \langle B_{\bar \gamma} (1 , \eta'(t_{-1})) \rangle $ are linearly independent if and only if $C_\gamma (1 , \eta'(t_{-1}))$ and $D_\gamma (1 , \eta'(t_{-1}))$ are linearly independent. In the sequence, we proceed by checking this condition. We observe that
$$C_\gamma \begin{pmatrix}
    1 \\ \eta'(t_{-1})
\end{pmatrix} = (L_{11}(t_{-1},t_0) + \eta'(t_{-1})) \begin{pmatrix}
    1 \\ L_{22} (t_0,t_1) 
\end{pmatrix}\, .$$
Moreover:
$$D_\gamma \begin{pmatrix}
    1 \\ \eta'(t_{-1})
\end{pmatrix} = -L_{12}^2(t_{-1},t_0) \begin{pmatrix}
    1 \\ L_{22} (t_0,t_1) 
\end{pmatrix} -L_{12}^2(t_0,t_1) (L_{11}(t_{-1},t_0) + \eta'(t_{-1})) \begin{pmatrix}
    0 \\ 1  
\end{pmatrix}\, .$$
By the twist condition, $L_{12} > 0$. Therefore, the two vectors above are lineraly independent if and only if there exists $t_{-1}$ such that $L_{11}(t_{-1},t_0) + \eta'(t_{-1}) \ne 0$, or equivalently:
$$L_{11}(t_{-1}, p_1 \circ T_\gamma (t_{-1}, \eta(t_{-1})) + \eta'(t_{-1}) \ne 0\, .$$
If this not the case, it means that the Lipschitz invariant graph $\eta$ satisfies a.e. the equation:
$$L_{11}(t_{-1},p_1 \circ T_\gamma (t_{-1}, \eta(t_{-1}))) + \eta'(t_{-1}) = 0$$
and consequently $\eta$ is at least $C^1$. However, in such a case, the $4$-periodic points $(t_i^{\pm},0)$ would be degenerate, since the stable and unstable directions should coincide, giving the desired contradiction. 

We have then showed that for a perturbation $\delta$ such that $\delta(t_0)=\delta'(t_0)=0$ and $\delta''(t_0)\neq 0$, there are $4$-periodic hyperbolic points within the Aubry-Mather set of rotation number $1/4$ with transverse heteroclinic or homoclinic intersections. In particular, the map $T_{\bar \gamma}$ associated to the perturbed domain does not have an invariant curve of rotation number $1/4$.
\end{proof}

The proof of Proposition \ref{open and dense no 1 4 inv curve} is concluded: it is sufficient to put together Lemma \ref{open no inv curve} and Proposition \ref{dense no inv curve}.

\section{Topologically and dynamically complex Birkhoff attractors}\label{birkhoff-section-complicated}

Birkhoff attractors described in Section \ref{NORMALLY} do not give an idea of their possible complexity, which in general occurs for mild dissipation (that is for $\lambda$ close to $1$). A sufficient condition in order to observe such a topologically and dynamically complex phenomena is that the conservative billiard map admits an instability region containing the zero section. In fact, by an adaptation of a result of Le Calvez (here Proposition \ref{lecalvez}), in such a case the Birkhoff attractor $\Lambda$ for the corresponding dissipative dynamics, with $\lambda$ close to 1, admits different upper and lower rotation numbers $\rho^- < \rho^+$ (defined in Section \ref{SUB 2.2}). This property has various consequences for $\Lambda$, all already observed in \cite{BFL} for Birkhoff billiards. In particular: 
\begin{itemize}
\item[$(i)$] $\Lambda$ is an \textit{indecomposable continuum}, that is it cannot be written as the union of two connected non-trivial sets (from a result by Charpentier, see \cite{Char34}). 
\item[$(ii)$] For every rational $p/q \in (\rho^- , \rho^+)$, there exists a periodic point in $\Lambda$ of the rotation number $\frac p q$ (which can be deduced from \cite{BG}).
\item[$(iii)$] If $x$ is a saddle periodic point of rotation number $\frac p q$, for $\frac p q\in(\rho^-,\rho^+)$, then its unstable manifold is contained in $\Lambda$ (see \cite[Proposition 14.3]{LeCal}).
\item[$(iv)$] The map $T_{\lambda}$ restricted to $\Lambda$ has positive topological entropy, as a consequence of the existence of a rotational horseshoe (see \cite[Theorem A]{Alejandro}).
\end{itemize}
We remind that an essential curve in $\mathbb{S} \times \mathbb{R}$ is a topological embedding of $\mathbb{S}$ which is not homotopic to a point. \\
\indent In the next proposition, we prove that, for a centrally symmetric table, any essential curve for the corresponding symplectic billiard dynamics passing through $\mathbb{S} \times \{0\}$ has necessarily rotation number $1/4$. 
\begin{proposition} \label{Hera}
Let $\Omega\subset \mathbb{R}^2$ be a $C^k$ strictly convex, centrally symmetric domain, $k\geq 2$. Denote by $T \colon \mathcal{P} \to \mathcal{P}$ the associated symplectic billiard map. Any essential invariant curve $\Gamma \subset \mathcal{P}$ passing through $\mathbb{S} \times \{0\}$ has rotation number $\rho(\Gamma) = 1/4$.
\end{proposition}

\begin{proof}
Recall that $\hat{T}: \hat{\mathcal{P}} \ni (t_0,t_1)\mapsto(t_1,t_2) \in \hat{\mathcal{P}}$ denotes the symplectic billiard map in the coordinates $(t_1,t_2)$. Let define
$$\hat{I}_2: \hat{\mathcal{P}} \ni (t_1,t_2) \mapsto (t_1,t_0^*) \in \hat{\mathcal{P}} $$
and
$$*: \hat{\mathcal{P}} \ni (t_1,t_2) \mapsto (t_1^*,t_2^*) \in \hat{\mathcal{P}}$$
Notice that $*$ is an involution, i.e., $*\circ*=Id$. Moreover, from the centrally symmetric assumption --see Figure \ref{I_2}-- we have that $\hat{T}$ and $*$ commute: this implies that also $\hat{I}_2$ and $*$ commute. Thus we obtain
$$\hat{T}\circ\hat{I}_2(t_{1},t_{2})=\hat{T}(t_{1},t_{0}^*)=(t_{0}^*,t_{-1})$$
and
$$\hat{I}_2\circ\hat{T}^{-1}(t_{1},t_{2})=\hat{I}_2(t_{0},t_{1})=(t_{0},t_{-1}^*)\, ,$$
so that, in particular: 
\begin{equation} \label{commute1}
\hat{T}\circ\hat{I}_2= * \circ \hat{I}_2 \circ \hat{T}^{-1} = \hat{I}_2\circ*\circ\hat{T}^{-1}\, .
\end{equation}


In the sequel, we simply denote by $\rho(\cdot)$ the rotation number with respect to $T$ and by $\rho(\cdot, T^2)$ the rotation number with respect to $T^2$. We divide the proof in two points. \\ 

\noindent $(i)$ We start by proving that $\rho(\Gamma)$ is necessarily rational. Let $\hat{\Gamma}\subset\hat{\mathcal{P}}$ be the image through $\phi: \mathcal{P} \to \hat{\mathcal{P}}$ --see formula (\ref{cambio})-- of such a curve $\Gamma$ and suppose --by contradiction-- that $\Gamma$ has irrational rotation number. First, we deduce, as a consequence of the centrally symmetric hypothesis, that also $*(\hat\Gamma)$ is invariant: indeed $\hat T \circ *=*\circ \hat T$ and, since $\hat T(\hat \Gamma)=\hat\Gamma$, the required invariance follows. Moreover, the two invariant curves $\hat\Gamma$ and $*(\hat\Gamma)$ must intersect, as we are going to show. With an abuse of notation, we denote by $*$ also the map corresponding to $*$ on $\mathcal{P}$, i.e.,
\[
*\colon (t,s)\in\mathcal{P}\to *(t,s)=(t^*,s)\in\mathcal{P}\, .
\]
Up to choose a good parametrization of the domain, we can assume that $t^*=t+\pi$ for every $t\in\mathbb{S}$. The curves $\Gamma,*(\Gamma)$ are both invariant, thus they are graphs. To fix the idea, we have
$$
\Gamma=\{(t,\eta(t)): \ t\in\mathbb{S}\}\quad \text{and}\quad *(\Gamma)=\{(t,\eta(t+\pi)): \ t\in\mathbb{S}\}\, .
$$
Since $\int_{\mathbb{S}}\eta(t)\, dt =\int_{\mathbb{S}}\eta(t+\pi)\, dt$, we deduce that $\Gamma\cap *(\Gamma)\neq \emptyset$.
Thus, the curve $*(\Gamma)$ has the same rotation number as $\Gamma$.
If, by contradiction, this rotation number is irrational, then $*(\Gamma)$ and $\Gamma$ 
must coincide (by e.g. \cite{Bang}[Section $4$]). Then, by (\ref{commute1}):
$$\hat{T}(\hat{I}_2(\hat\Gamma))=\hat{I}_2\circ*\circ\hat{T}^{-1}(\hat\Gamma)=\hat{I}_2(\hat\Gamma)\, .$$
Observe that the map corresponding to $\hat{I_2}$ in $\mathcal{P}$ is the map $I_2\colon (t,s)\in\mathcal{P}\to (t,-s)\in\mathcal{P}$. The image corresponding to the curve $\hat{I}_2(\hat\Gamma)$ is then $I_2(\Gamma)$.
The previous equality means that also $I_2(\Gamma)$ is an invariant essential curve and therefore, since $\Gamma$ is a graph,
$$\Gamma \cap \left(\mathbb{S} \times \{0\}\right) \ne \emptyset \Leftrightarrow \Gamma \cap I_2(\Gamma) \ne \emptyset\, .$$
Since $\Gamma$ and $I_2(\Gamma)$ intersect, these two curves have the same rotation number, i.e., $\rho(\Gamma) = \rho(I_2(\Gamma))$. Again, since the rotation number is supposed irrational, we deduce that $\Gamma=I_2(\Gamma)=\mathbb{S}\times\{0\}$. This contradicts the fact that the zero section contains $4$ periodic points and, thus, that its rotation number should be $1/4 \in \mathbb{Q}$.\\ \\ 
\noindent $(ii)$ We finally prove that $\rho(\Gamma)$ is necessarily $1/4$. We start by recalling that $\rho(\Gamma) = \rho(*(\Gamma))$, since they are both invariant curves and the intersection of $\Gamma$ and $*(\Gamma)$ is not empty, as explained above. Denote by $p/q$ the rotation number of $\Gamma$. If $\Gamma \cap *(\Gamma) \subset \mathbb{S} \times \{ 0 \}$, then $\Gamma \cap *(\Gamma)$ in an invariant set necessarily given by $4$-periodic points. Hence $\rho(\Gamma) = 1/4$. We then suppose --on the contrary-- that $\Gamma \cap *(\Gamma) \not \subset \mathbb{S} \times \{ 0 \}$. To conclude, we need to consider also the essential curve $I_2(\Gamma)$, which is $T^2$-invariant. In fact, applying twice \eqref{commute1} and since $\hat T$ and $*$ commute, we get:  
$$ \hat{T}^2 \circ \hat{I}_2(\hat{\Gamma}) = \hat{T} \circ * \circ \hat{I}_2(\hat{\Gamma}) = * \circ \hat{T} \circ \hat{I}_2(\hat{\Gamma}) = * \circ * \circ \hat{I}_2(\hat{\Gamma}) = \hat{I}_2(\hat{\Gamma})\, .$$
From now on, with abuse of notation, we continue to indicate with $T$ its lift to $\mathbb{R}^2$. By Proposition \ref{T dopo T twist}, $T^2$ is a twist map so we can consider the rotation number of the essential $T^2$-invariant curves $\Gamma$, $*(\Gamma)$ and $I_2(\Gamma)$:
$$\frac{2p}{q} = 2 \rho(\Gamma) = \rho(\Gamma, T^2) = \rho(*(\Gamma), T^2) = \rho(I_2(\Gamma), T^2)\, .$$
In the last two equalities, we have respectively used the facts that the rotation numbers are the same because both $\Gamma \cap *(\Gamma)$ and $\Gamma \cap I_2(\Gamma)$ are nonempty. \\
\noindent By using classical Aubry-Mather theory (see \cite{Bang} and also \cite[Subsection 3.4]{MCA}), in the sequel we will prove that
$\Gamma \cap *(\Gamma) \cap I_2(\Gamma) \ne \emptyset$ since this intersection necessarily contains every $2p/q$-periodic point for $T^2$. \\ 
Recall that $T^2$ is still a twist map and that we are denoting by $T^2$ also a lift of $T^2$. Since $\Gamma$, $*(\Gamma)$ and $I_2(\Gamma)$ are all $T^2$-invariant graphs, they consist of points of $T^2$-minimizing orbits (see \cite[Proposition 2.8]{Mather}, \cite{Bang} or \cite[Subsection 5.2]{MCA}). Equivalently, $\Gamma, *(\Gamma), I_2(\Gamma)$ are contained in
$$\mathcal{M}_{\frac{2p}{q}}(T^2) := \{(t,s) \in \mathcal{P} \text{ having a $T^2$-minimizing orbit of rotation number $2p/q$}\}\, .$$
We recall that set $\mathcal{M}_{\frac{2p}{q}}(T^2)$ is the disjoint union of $3$ invariant sets:
$$\mathcal{M}_{\frac{2p}{q}}(T^2) = \mathcal{M}_{\frac{2p}{q}}^{per}(T^2) \sqcup \mathcal{M}_{\frac{2p}{q}}^{+}(T^2) \sqcup \mathcal{M}_{\frac{2p}{q}}^{-}(T^2)\, ,$$
where, denoting by $p_1$ the projection on the first coordinate,
$$\mathcal{M}_{\frac{2p}{q}}^{per}(T^2) = \{(t,s) \in \mathcal{M}_{\frac{2p}{q}}(T^2)| \ p_1 \circ T^{2q}(t,s) = p_1(t,s) + 2p \}\, ,$$ 
$$\mathcal{M}_{\frac{2p}{q}}^{+}(T^2) = \{(t,s) \in \mathcal{M}_{\frac{2p}{q}}(T^2)| \ p_1 \circ T^{2q}(t,s) > p_1(t,s) + 2p \}\, ,$$ 
and
$$\mathcal{M}_{\frac{2p}{q}}^{-}(T^2) = \{(t,s) \in \mathcal{M}_{\frac{2p}{q}}(T^2)| \ p_1 \circ T^{2q}(t,s) < p_1(t,s) + 2p\}\, .$$
\noindent Since both $\mathcal{M}_{\frac{2p}{q}}^{per}(T^2)\sqcup \mathcal{M}_{\frac{2p}{q}}^{+}(T^2)$ and $\mathcal{M}_{\frac{2p}{q}}^{per}(T^2)\sqcup \mathcal{M}_{\frac{2p}{q}}^{-}(T^2)$ are well-ordered sets, we have that the fiber of every $P\in \mathcal{M}_{\frac{2p}{q}}^{per}(T^2)$ intersects the set $\mathcal{M}_{\frac{2p}{q}}(T^2)$ uniquely in $P$. This means that:
\begin{equation}\label{uguaglianza insiemi}p_1^{-1}\left(p_1\left(\mathcal{M}_{\frac{2p}{q}}^{per}(T^2)\right) \right) \bigcap \mathcal{M}_{\frac{2p}{q}}(T^2) = \mathcal{M}_{\frac{2p}{q}}^{per}(T^2)\, ,
\end{equation}
see \cite[Section 5]{Bang}. We deduce that the set $\mathcal{M}_{\frac{2p}{q}}^{per}(T^2)$ is necessarily contained in every essential $T^2$-invariant curve of rotation number $2p/q$. In fact, let $\Gamma$ be an essential $T^2$-invariant curve of rotation number $\frac{2p}{q}$: then, it is contained in $\mathcal{M}_{\frac{2p}{q}}(T^2)$ and it is a graph $\{ (t,\eta(t)):\ t\in\mathbb{R}\}$. Let $P\in \mathcal{M}_{\frac{2p}{q}}^{per}(T^2)$ and let $p_1(P)$ be its projection on the first coordinate. From \eqref{uguaglianza insiemi}, we deduce that the point $(p_1(P),\eta(p_1(P))\in\Gamma$ must be the point $P$, i.e., $P\in\Gamma$.
In particular:


$$\emptyset \ne \mathcal{M}_{\frac{2p}{q}}^{per}(T^2) \subseteq \Gamma \cap *(\Gamma) \cap I_2(\Gamma)\, .$$
Since $\Gamma\cap I_2(\Gamma)\subset\mathbb{S}\times\{0\}$, $\mathcal{M}_{\frac{2p}{q}}^{per}(T^2)$ is given by periodic points for $T^2$ contained in $\mathbb{S}\times\{0\}$. Given $P\in\mathcal{M}_{\frac{2p}{q}}^{per}(T^2)$, the $T^2$-orbit of $T(P)$ is clearly contained in $\Gamma$, since $\Gamma$ is actually $T$-invariant. Consequently, since $T(P)$ is periodic for $T^2$, its orbit belongs to $\mathcal{M}_{\frac{2p}{q}}^{per}(T^2)$ and therefore it is contained in $\mathbb{S}\times\{0\}$.  This means that the $T$-orbit of $P$ is entirely contained in the zero section, and therefore $\rho(\Gamma)=1/4$.

 \end{proof}

In the sequel, the union of all $T$-invariant essential curves in $\mathcal{P}$ will be indicated by $\mathcal{V}(T)$. The next definition can be formulated for a general twist map (see e.g. Definition 6.9 in \cite{BFL}). 
\begin{definition}
An instability region for $T: \mathcal{P} \to \mathcal{P}$ is an open bounded connected component of $\mathcal{P} \setminus \mathcal{V}(T)$ that contains in its interior an essential curve. 
\end{definition}
\noindent We conclude this preamble by recalling a result in \cite{BFL} (see Proposition 6.10 and Corollary 6.12) whose proof adapts a former argument due to Le Calvez (see \cite[Section 8]{LeCal}). To do that, we remind the notion of twist map with respect to $\beta \in (0,\frac{\pi}{2})$ (given for general twist maps in \cite[Section 1.2]{Herman}, see also \cite{LeCalEns} or \cite[Section 6.1]{BFL}). 
\begin{definition} \label{lecalvez docet} Let $U \subset \mathcal{P}$ be open.
We say that $T: U \to T(U)$ is a positive (resp. negative) twist map on $U$ with respect to $\beta \in (0,\frac{\pi}{2})$ if it is a $C^1$ diffeomorphism onto its image and, for every $(t,s) \in U$, the angle formed by the unitary vertical vector $(0,1) \in T_{(t,s)} \mathcal{P}$ and $DT(t,s)(0,1)$ is in $(\beta - \pi,- \beta)$ (resp. $(\beta, \pi - \beta)$), where we have fixed at every tangent plane $T_{(t,s)}\mathcal{P}$ the counter-clockwise orientation. 
\end{definition}

\begin{figure}[ht]
    \centering
    \includegraphics[width=0.3\linewidth]{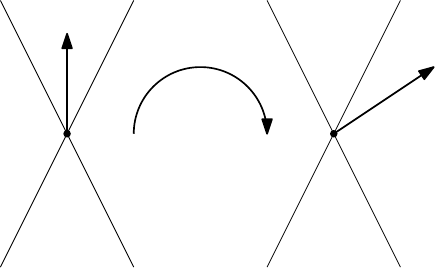}
    \caption{Definition of twist map.}
    \label{fig:ConeTwist}
\end{figure}

\noindent We refer to Figure \ref{fig:ConeTwist}.

\begin{remark}
We notice that, whenever the boundary of the symplectic billiard table is $C^2$ strictly convex, $T$ is a negative twist map according to the previous Definition \ref{lecalvez docet}. In fact, let $(t_1, s_1) \in \mathcal{P}$. The image of the vertical direction $(0,1)$ by the differential of $T$ --see formula \eqref{differential}-- is
$$DT(t_1,s_1)(0,1) = -\frac{1}{L_{12}(t_1,t_2)}\begin{pmatrix}
1 \\
L_{22} (t_1,t_2)
\end{pmatrix}
$$
and, independently from the point $(t_1, s_1) \in \mathcal{P}$, we have that $L_{22} (t_1,t_2)$ is uniformly bounded.   \end{remark}
\begin{proposition} \label{lecalvez} Let $\Omega$ be a $C^k$ strictly convex domain, $k \ge 2$. Suppose that there exists $$\mathcal{I} := \{ (t,s): \ \phi^-(t) < s < \phi^+(t) \} \subseteq \mathcal{P}, \qquad \phi^- < \phi^+$$ 
where $\phi^{\pm}: \mathbb{S} \to \mathbb{R}$ are continuous maps such that:
\begin{itemize}
\item[$(a)$] $T: \mathcal{I} \to \mathcal{I}$ is a positive twist map with respect to $\beta \in (0,\frac{\pi}{2})$.
\item[$(b)$] $\mathcal{I}$ is an instability region for $T$ containing $\mathbb{S} \times \{0\}$. 
\end{itemize}
Then there exist $\lambda_0 \in (0,1)$ such that, for any $\lambda \in [\lambda_0,1)$, the Birkhoff attractor $\Lambda$ of the corresponding dissipative symplectic billiard map $T_{\lambda}$ has $\rho^+-\rho^->0$, with $\frac{1}{4} \in (\rho^-,\rho^+) \mod \mathbb{Z}$.
\end{proposition}
\noindent The rest of the section is devoted to proving that hypothesis of Proposition \ref{lecalvez} are satisfied in the following two cases.
\begin{itemize}
\item[$(a)$] For an open and dense set of strongly convex, centrally symmetric billiards tables. Indeed, using Proposition \ref{open and dense no 1 4 inv curve}, we will see that there exists an open and dense set of tables exhibiting a region of instability containing the zero section, see Theorem \ref{N PER}.
\item[$(b)$] For symplectic billiard maps when the strictly convex billiard table has a point with zero curvature. This result --here Proposition \ref{mather docet}-- is a straightforward consequence of Mather's theorem on non-existence of caustics which holds also for symplectic billiards (as proved in \cite[Theorem 2]{AT}). 
\end{itemize}

\indent We proceed with precise statements and proofs of $(a)$ and $(b)$, which are Theorem \ref{N PER} and Proposition \ref{mather docet}, respectively. In particular, Theorem \ref{N PER} immediately follows from Proposition \ref{Hera}.
\begin{theorem} \label{N PER}
The set of strongly convex, centrally symmetric billiard tables whose associated
symplectic billiard map does have an instability region containing $\mathbb{S}\times\lbrace0\rbrace$ is open and dense among $\mathcal{B}_2$.
\end{theorem}
\begin{proof}
According to Proposition \ref{Hera}, any essential invariant curve $\Gamma \subset \mathcal{P}$ passing
through $\mathbb{S} \times \{0\}$ has rotation number $\rho(\Gamma) = 1/4$. To conclude is then sufficient to apply Proposition \ref{open and dense no 1 4 inv curve}, assuring that the set of strongly convex, centrally symmetric billiard tables whose associated symplectic billiard map does not have an invariant curve of rotation number $1/ 4$, is open and dense among $\mathcal{B}_2$.
\end{proof}

\begin{proposition} \label{mather docet}
If the curvature of the boundary of a $C^2$ strictly convex billiard table vanishes at some point, then the whole associated $\mathcal{P}$ is an instability region. 
\end{proposition}
\begin{proof}
According to Theorem 2 in \cite{AT}, when the curvature of the boundary of the billiard table vanishes at some point, then the associate symplectic billiard map $T$ has no caustics. This means that the whole phase space $\mathcal{P}$ is an instability region.
\end{proof}

\section{Numerical simulations} \label{NUM SIM}
In this section, we present some numerical simulations to illustrate the results discussed in the previous sections. With the aid of \textit{Mathematica}, we compute the billiard map $T_\lambda$ for specific domains, both centrally symmetric and non-symmetric, and plot some orbits in the corresponding phase space. 

To this end, we represent the domain using the angle $\theta$, defined as the angle between the tangent vector and a fixed reference direction, along with the support function as defined in Section 4. We then choose to represent the phase space using the coordinates $\theta$ and $\psi = \theta_2 - \theta_1$, that is, the angular difference between two consecutive points of the orbit (up to consider suitable lifts of angles). As a result, the phase space reduces to $\mathbb{S} \times [0, \pi]$.

We then plot several orbits $T^n_\lambda(\theta_0, \psi_0)$, choosing random initial values for $\theta_0$ and $\psi_0$. Each orbit is assigned a different color according to its initial value $\psi_0$, with cooler colors corresponding to values of $\psi_0$ close to 0, and progressively warmer colors assigned as $\psi_0$ increases.

The obtained simulations let appear the attractor (not necessarily the Birkhoff attractor). To gain some insight into the structure of the attractor, we display only the points of the orbits for $n > n_0$, where $n_0$ is chosen appropriately.

\begin{figure}[ht]
    \centering
    \includegraphics[width=0.5\linewidth]{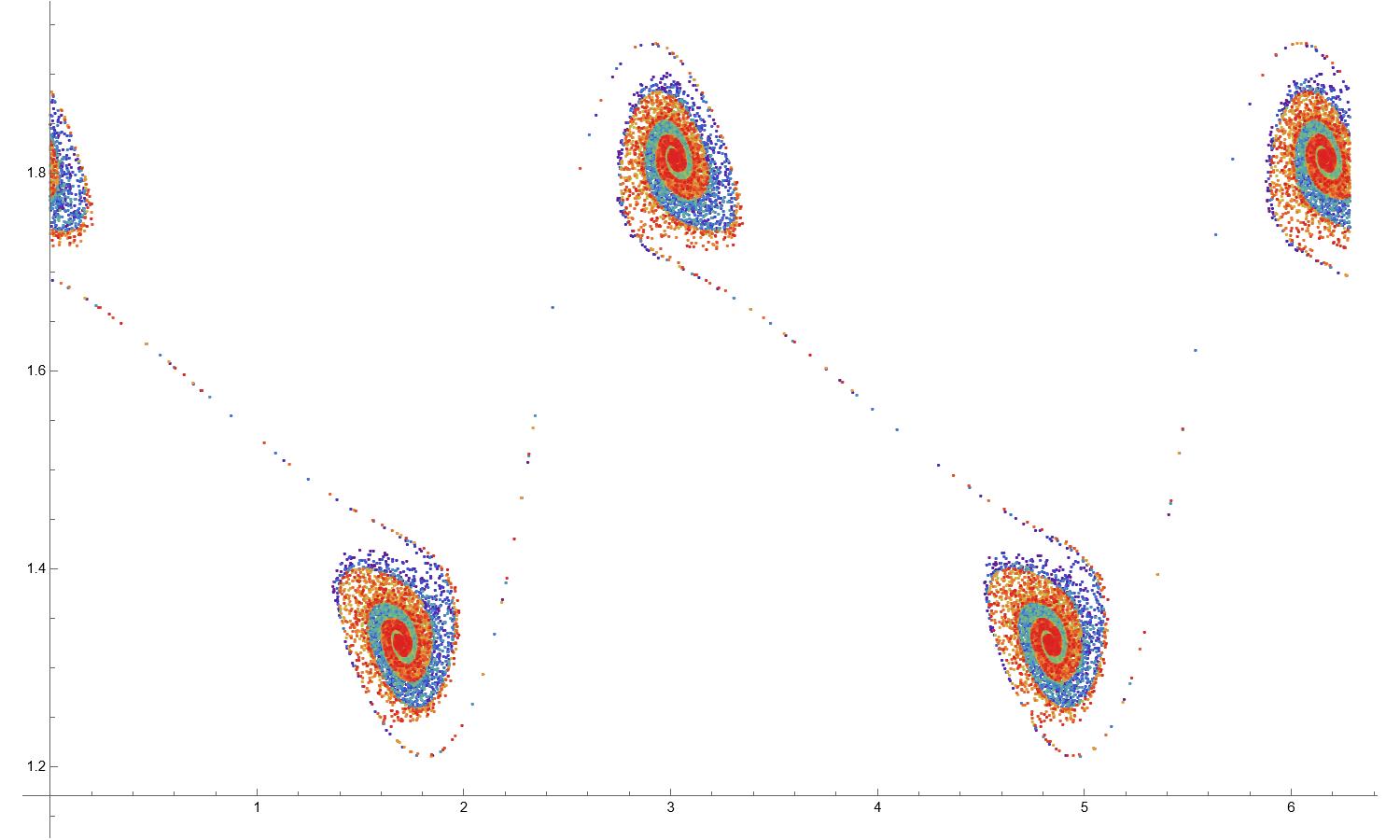}
    \caption{ $p(\theta)=1+\frac{\sin{2\theta}}{8}$, $\lambda=0.9$, $n_0=30$.}
    \label{num1}
\end{figure}

 In Figure \ref{num1}, we consider a non-symmetric billiard table. From the simulation, it seems that the attractor and the Birkhoff attractor coincide, even if the latter is not a graph over $\mathbb{S}$. A clearly visible 4-periodic orbit can be identified, around which the attractor wraps.
 
\begin{figure}[ht]
    \centering
\includegraphics[width=0.5\linewidth]{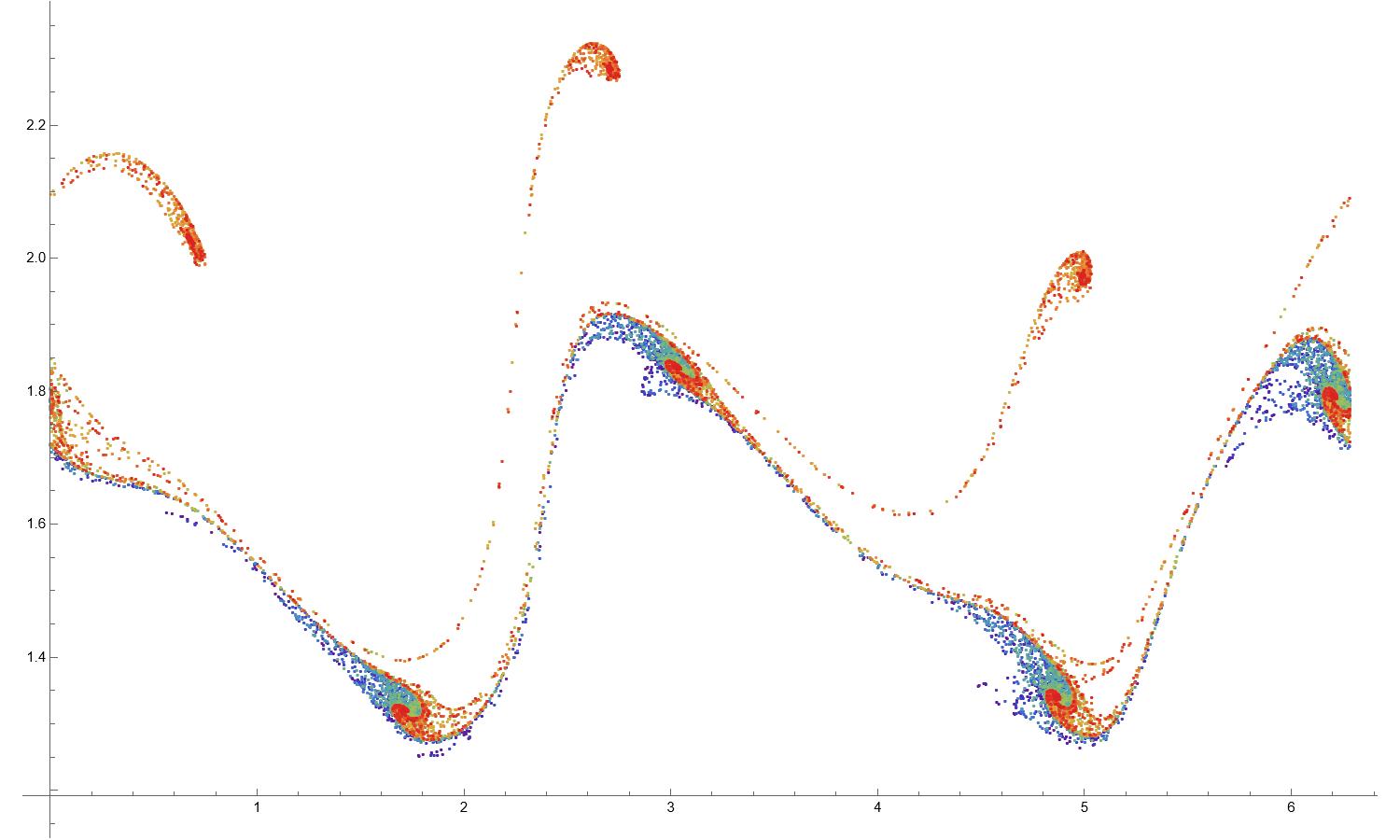}
    \caption{ $p(\theta)=1+\frac{\sin{2\theta}}{8}+\frac{\cos{3\theta}}{27}$, $\lambda=0.71$, $n_0=10$.}
    \label{num2}
\end{figure}

\begin{figure}[H]
    \centering
\includegraphics[width=0.5\linewidth]{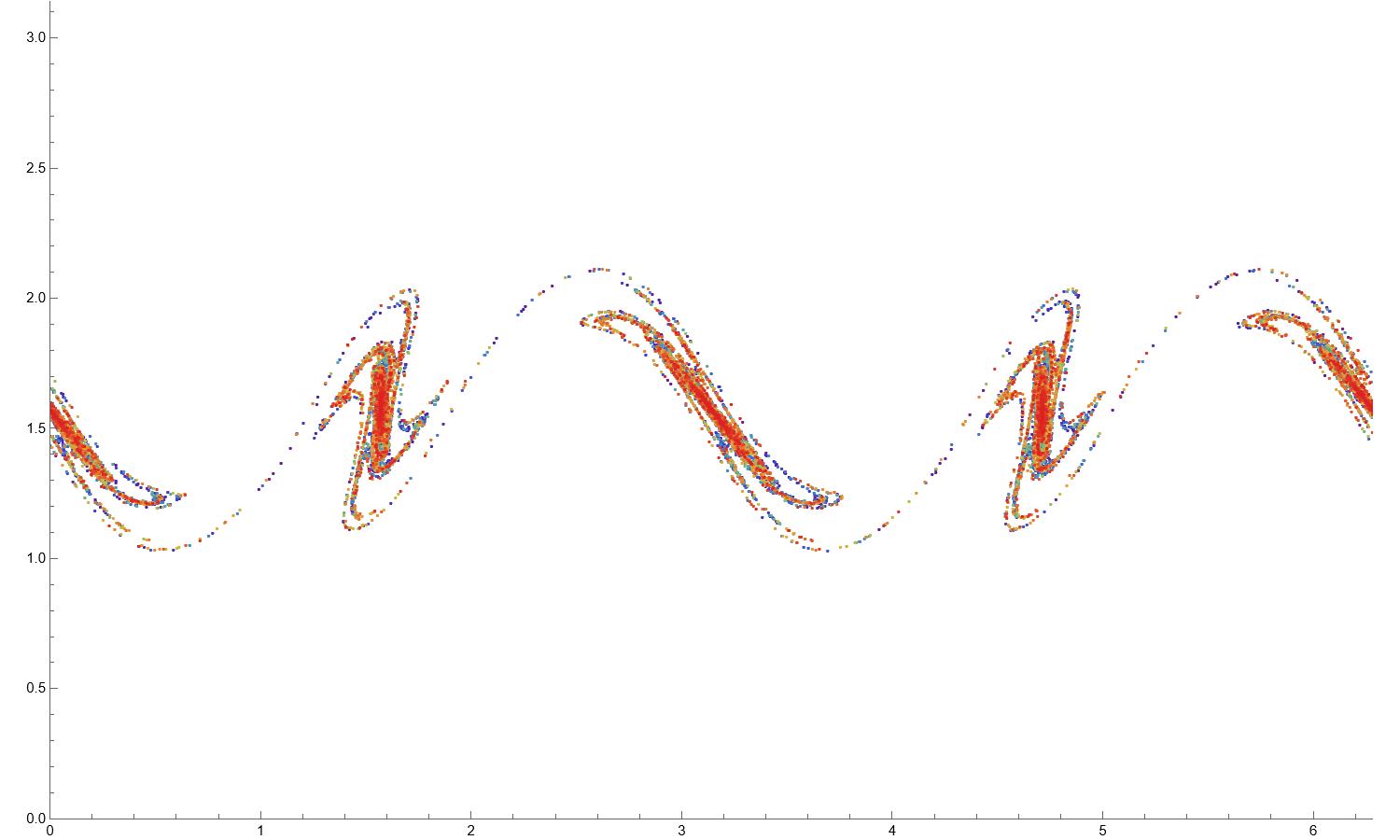}
    \caption{ $r(\theta)=1-\frac{\cos{2\theta}}{5}$, $\lambda=0.71$, $n_0=10$.}
    \label{num3}
\end{figure}

Figure \ref{num2}, again for a non-symmetric table, by contrast, seems to give an example in which the attractor strictly contains the Birkhoff attractor. There are extra components in the attractor, which, we guess, are due to the presence of a 3-periodic orbit. \\
\indent 
In Figure \ref{num3}, we consider the case of a billiard table that is centrally symmetric and has two points of zero curvature. In this example, we use polar coordinates to represent the table.
It is well known from the previous section, that in this case, the entire phase space forms a region of instability, and that for dissipation values close to 1, the attractor becomes topologically complex. Notably, there are blue-colored points in the upper part and red-colored points in the lower part of the figure. This indicates that the attractor is highly intricate and entangled.

\bibliographystyle{plain}
\bibliography{biblio.bib}

\begin{thebibliography}{10}

\bibitem{AT}
P.~Albers and S.~Tabachnikov.
\newblock Introducing symplectic billiards.
\newblock {\em Adv. Math.}, 333:822--867, 2018.

\bibitem{AA}
S.~Allais and M.-C. Arnaud.
\newblock The dynamics of conformal {H}amiltonian flows: dissipativity and
  conservativity.
\newblock {\em Rev. Mat. Iberoam.}, 40(3):987--1021, 2024.

\bibitem{MCA}
M.-C. Arnaud.
\newblock Hyperbolicity for conservative twist maps of the 2-dimensional
  annulus.
\newblock {\em Publ. Mat. Urug.}, 16:1--39, 2016.

\bibitem{AF}
M.-C. Arnaud and J.~Fejoz.
\newblock Invariant submanifolds of conformal symplectic dynamics.
\newblock {\em J. \'Ec. polytech. Math.}, 11:159--185, 2024.

\bibitem{AFR}
M.-C. Arnaud, A.~Florio, and V.~Roos.
\newblock Vanishing asymptotic {M}aslov index for conformally symplectic flows.
\newblock {\em Ann. H. Lebesgue}, 7:307--355, 2024.

\bibitem{AHV}
M.-C. Arnaud, V.~Humilière, and C.~Viterbo.
\newblock {H}igher {D}imensional {B}irkhoff attractors, 2024.
\newblock arXiv:2404.00804.

\bibitem{VitPlusPlus}
T.~Asano, S.~Guillermou, Y.~Ike, and C.~Viterbo.
\newblock Regular {L}agrangians are smooth {L}agrangians, 2024.
\newblock arXiv:2407.00395, to appear in {J}. {M}ath. {S}oc. {J}apan.

\bibitem{Bang}
V.~Bangert.
\newblock Mather sets for twist maps and geodesics on tori.
\newblock In {\em Dynamics reported, {V}ol.\ 1}, volume~1 of {\em Dynam.
  Report. Ser. Dynam. Systems Appl.}, pages 1--56. Wiley, Chichester, 1988.

\bibitem{BaBe}
L.~Baracco and O.~Bernardi.
\newblock Totally integrable symplectic billiards are ellipses.
\newblock {\em Adv. Math.}, 454:Paper No. 109873, 17, 2024.

\bibitem{BBN2}
L.~Baracco, O.~Bernardi, and A.~Nardi.
\newblock {B}ialy-{M}ironov type rigidity for centrally symmetric symplectic
  billiards.
\newblock {\em Nonlinearity}, 37(12):Paper No. 125025, 12, 2024.

\bibitem{BBN1}
L.~Baracco, O.~Bernardi, and A.~Nardi.
\newblock Higher order terms of {M}ather's {$\beta $}-function for symplectic
  and outer billiards.
\newblock {\em J. Math. Anal. Appl.}, 537(2):Paper No. 128353, 20, 2024.

\bibitem{BBN3}
L.~Baracco, O.~Bernardi, and A.~Nardi.
\newblock Area spectral rigidity for axially symmetric and {Radon} domains.
\newblock Preprint, {arXiv}:2410.12644 [math.{DS}] (2025), 2025.

\bibitem{BG}
M.~Barge and R.~M. Gillette.
\newblock Rotation and periodicity in plane separating continua.
\newblock {\em Ergodic Theory Dynam. Systems}, 11(4):619--631, 1991.

\bibitem{BergerBounemoura}
P.~Berger and A.~Bounemoura.
\newblock A geometrical proof of the persistence of normally hyperbolic
  submanifolds.
\newblock {\em Dyn. Syst.}, 28(4):567--581, 2013.

\bibitem{BFL}
O.~Bernardi, A.~Florio, and M.~Leguil.
\newblock {B}irkhoff attractors of dissipative billiards.
\newblock {\em Ergodic Theory and Dynamical Systems}, 2024.
\newblock https://doi.org/10.1017/etds.2024.68.

\bibitem{B22}
G.~D. Birkhoff.
\newblock Surface transformations and their dynamical applications.
\newblock {\em Acta Math.}, 43(1):1--119, 1922.

\bibitem{Birkhoff}
G.~D. Birkhoff.
\newblock Sur quelques courbes ferm\'ees remarquables.
\newblock {\em Bull. Soc. Math. France}, 60:1--26, 1932.

\bibitem{CCDLL}
R.~Calleja, A.~Celletti, and R.~de~la Llave.
\newblock K{AM} theory for some dissipative systems.
\newblock In {\em New frontiers of celestial mechanics---theory and
  applications}, volume 399 of {\em Springer Proc. Math. Stat.}, pages 81--122.
  Springer, Cham, [2022] \copyright 2022.

\bibitem{PassTalPlus}
M.~Capi\`nski, M.~Gröger, A.~Passeggi, and F.~A. Tal.
\newblock Conditions {I}mplying {A}nnular {C}haos: {Q}uantitative results and
  {C}omputer {A}ssisted {P}roofs, 2025.
\newblock arXiv:2506.06608.

\bibitem{Char34}
M.~Charpentier.
\newblock Sur quelques propriétés des courbes de m. birkhoff.
\newblock {\em Bulletin de la Société Mathématique de France}, 62:193--224,
  1934.

\bibitem{Zavi3}
Q.~Chen, A.~Fathi, M.~Zavidovique, and J.~Zhang.
\newblock Convergence of the solutions of the nonlinear discounted
  {H}amilton-{J}acobi equation: the central role of {M}ather measures.
\newblock {\em J. Math. Pures Appl. (9)}, 181:22--57, 2024.

\bibitem{Cro}
S.~Crovisier.
\newblock Langues d'{A}rnold{} g\'en\'eralis\'ees des applications de l'anneau
  d\'eviant la verticale.
\newblock {\em C. R. Math. Acad. Sci. Paris}, 334(1):47--52, 2002.

\bibitem{CroPot}
S.~Crovisier and R.~Potrie.
\newblock Introduction to partially hyperbolic dynamics.
\newblock ICTP, Trieste (2015).

\bibitem{Zavi}
A.~Davini and M.~Zavidovique.
\newblock Convergence of the solutions of discounted {H}amilton-{J}acobi
  systems.
\newblock {\em Adv. Calc. Var.}, 14(2):193--206, 2021.

\bibitem{Day}
M.~M. Day.
\newblock Polygons circumscribed about closed convex curves.
\newblock {\em Trans. Amer. Math. Soc.}, 62:315--319, 1947.

\bibitem{FSV}
C.~Fierobe, A.~Sorrentino, and A.~Vig.
\newblock Deformational spectral rigidity of axially-symmetric symplectic
  billiards.
\newblock Preprint, {arXiv}:2410.13777 [math.{DS}] (2024), 2024.

\bibitem{Flan68}
H.~Flanders.
\newblock A proof of minkowski's inequality for convex curves.
\newblock {\em The American Mathematical Monthly}, 75(6):581--593, 1968.

\bibitem{GDLLS}
M.~Gidea, R.~de~la Llave, and T.~M-Seara.
\newblock Geometric, topological and dynamical properties of conformally
  symplectic systems, normally hyperbolic invariant manifolds, and scattering
  maps, 2025.
\newblock arXiv:2508.14794.

\bibitem{Herman}
M.~Herman.
\newblock {\em Sur les courbes invariantes par les diff\'eomorphismes de
  l'anneau. {V}ol. 1}, volume 103-104 of {\em Ast\'erisque}.
\newblock Soci\'et\'e{} Math\'ematique de France, Paris, 1983.
\newblock With an appendix by Albert Fathi, With an English summary.

\bibitem{HPS}
M.~W. Hirsch, C.~C. Pugh, and M.~Shub.
\newblock Invariant manifolds.
\newblock {\em Bull. Amer. Math. Soc.}, 76:1015--1019, 1970.

\bibitem{Koropecki}
A.~Koropecki.
\newblock Realizing rotation numbers on annular continua.
\newblock {\em Math. Z.}, 285(1-2):549--564, 2017.

\bibitem{LeCal}
P.~Le~Calvez.
\newblock Propri\'et\'es des attracteurs de {B}irkhoff.
\newblock {\em Ergodic Theory Dynam. Systems}, 8(2):241--310, 1988.

\bibitem{LeCalEns}
P.~Le~Calvez.
\newblock Étude topologique des applications déviant la verticale. / patrice
  le calvez.
\newblock {\em Ensaios Matemáticos}, 2:7--102, 1990.

\bibitem{MPS}
R.~Markarian, E.~J. Pujals, and M.~Sambarino.
\newblock Pinball billiards with dominated splitting.
\newblock {\em Ergodic Theory Dynam. Systems}, 30(6):1757--1786, 2010.

\bibitem{MS}
S.~Mar\`o and A.~Sorrentino.
\newblock Aubry-{M}ather theory for conformally symplectic systems.
\newblock {\em Comm. Math. Phys.}, 354(2):775--808, 2017.

\bibitem{MaSw}
H.~Martini and K.~J. Swanepoel.
\newblock Equiframed curves---a generalization of {R}adon curves.
\newblock {\em Monatsh. Math.}, 141(4):301--314, 2004.

\bibitem{Mather}
J.~N. Mather.
\newblock Variational construction of orbits of twist diffeomorphisms.
\newblock {\em J. Amer. Math. Soc.}, 4(2):207--263, 1991.

\bibitem{Alejandro}
A.~Passeggi, R.~Potrie, and M.~Sambarino.
\newblock Rotation intervals and entropy on attracting annular continua.
\newblock {\em Geom. Topol.}, 22(4):2145--2186, 2018.

\bibitem{PassTal}
A.~Passeggi and F.~A. Tal.
\newblock Conditions implying annular chaos, 2023.
\newblock arXiv:2305.02963.

\bibitem{Sambarino}
M.~Sambarino.
\newblock A (short) survey on dominated splittings.
\newblock In {\em Mathematical {C}ongress of the {A}mericas}, volume 656 of
  {\em Contemp. Math.}, pages 149--183. Amer. Math. Soc., Providence, RI, 2016.

\bibitem{Tab}
S.~Tabachnikov.
\newblock Outer billiards.
\newblock {\em Uspekhi Mat. Nauk}, 48(6(294)):75--102, 1993.

\bibitem{Tsod}
D.~Tsodikovich.
\newblock Local rigidity for symplectic billiards.
\newblock {\em J. Geom. Anal.}, 35(10):26, 2025.
\newblock Id/No 306.

\bibitem{Vit22}
C.~Viterbo.
\newblock On the supports in the {H}umilière completion and
  $\gamma$-coisotropic sets, 2022.
\newblock arXiv:2204.04133.

\bibitem{Zavi2}
M.~Zavidovique.
\newblock Convergence of solutions for some degenerate discounted
  {H}amilton-{J}acobi equations.
\newblock {\em Anal. PDE}, 15(5):1287--1311, 2022.

\end{thebibliography}
\end{document}